\documentclass[11pt]{article}
\usepackage{amsmath, graphicx, amsfonts,amssymb}
\usepackage{amsfonts}
\usepackage{mathrsfs}
\usepackage{amsthm}
\usepackage{cite}
\usepackage[pdfstartview=FitH,
            CJKbookmarks=true,
            bookmarksnumbered=true,
            bookmarksopen=true,
            colorlinks=true,
            linkcolor=blue,
            anchorcolor=blue,
            citecolor=red,
            urlcolor=blue
            ]{hyperref}

\numberwithin{equation}{section}
\newtheorem{theorem}{Theorem}[section]
\newtheorem{definition}[theorem]{Definition}
\newtheorem{lemma}[theorem]{Lemma}
\newtheorem{corollary}[theorem]{Corollary}

\newtheorem{proposition}[theorem]{Proposition}

\def\bl{\begin{lemma}}
\def\el{\end{lemma}}
\def\bc{\begin{corollary}}
\def\ec{\end{corollary}}
\def\bt{\begin{theorem}}
\def\et{\end{theorem}}
\def\Cyk{\mathcal{C}_{e}^n}

\def\bp{\begin{proposition}}
\def\ep{\end{proposition}}
\def\be{\begin{equation}}
\def\ee{\end{equation}}
\def\bd{\begin{definition}}
\def\ed{\end{definition}}
\def\plp{\Lambda_p^{\circ}K}
\def\pli{\Lambda_{\infty}^{\circ}K}

\def\ksp{K^{\diamond}}

\newcommand{\R}{\mathbb{R}}
\newcommand{\Rn}{\mathbb{R}^n}
\def\sphere{S^{n-1}}
\def\ball{B^n}

\def\cKon{\mathscr{K}_{o}^n}
\def\cKe{\mathscr{K}_{e}^n}

\def\polar{K^\circ}

\def\cS{\mathcal{S}}

\def\bbE{\mathbb{E}}

\title{On the  sine
polarity and the $L_p$-sine Blaschke-Santal\'{o} inequality
 \footnote{Keywords: Blaschke-Santal\'{o} inequality, $L_p$-sine transform, $L_p$-sine centroid body, sine polar body.}
 }

 \author{Qingzhong Huang, Ai-Jun Li \footnote{Corresponding author: Ai-Jun Li}, Dongmeng Xi and  Deping Ye}
\begin{document}

\date{}
\maketitle

\begin{abstract} This paper is dedicated to study the sine version of polar bodies and establish the $L_p$-sine Blaschke-Santal\'{o} inequality
for the $L_p$-sine centroid body. 

The $L_p$-sine centroid body
$\Lambda_p K$ for a star body  $K\subset\Rn$ is a convex body based
on the $L_p$-sine transform, and its associated
Blaschke-Santal\'{o} inequality provides an upper bound for the
volume of $\Lambda_p^{\circ}K$, the polar body of $\Lambda_p K$, in
terms of the volume of $K$. Thus, this inequality can be viewed as
the ``sine cousin" of the $L_p$ Blaschke-Santal\'{o} inequality
established by Lutwak and Zhang. As $p\rightarrow \infty$, the limit of $\Lambda_p^{\circ} K$ becomes
the sine polar body $K^{\diamond}$ and hence the $L_p$-sine
Blaschke-Santal\'{o} inequality reduces to the sine
Blaschke-Santal\'{o} inequality for the sine polar body.  The sine
polarity naturally leads to a new class of convex bodies
$\mathcal{C}_{e}^n$, which consists of all origin-symmetric convex
bodies generated by the intersection of origin-symmetric closed
solid cylinders. Many notions in $\mathcal{C}_{e}^n$ are developed,
including the cylindrical support function, the supporting cylinder,
the cylindrical Gauss image, and the cylindrical hull. Based on
these newly introduced notions, the equality conditions of the sine
Blaschke-Santal\'{o} inequality are settled.

 \vskip 2mm \noindent  2020 Mathematics Subject Classification: 52A20, 52A30, 52A40, 94A15.
\end{abstract}

\section{Introduction}

Geometric inequalities, such as isoperimetric or affine
isoperimetric inequalities, are the central objects of interest in
convex geometry. These inequalities aim to estimate geometric (or
affine) invariants from above and/or below in terms of the volume.
One of the most important affine isoperimetric inequalities is the
celebrated Blaschke-Santal\'{o} inequality. It asserts that,  in the $n$-dimensional Euclidean space $\Rn$,
\begin{equation}\label{BS-polar-1} V(K)V(\polar)\leq
\omega_n^2\end{equation} holds for any origin-symmetric convex body
$K$ or more generally for any convex body $K$ with its centroid at the origin,  with equality if and
only if $K$ is an origin-symmetric ellipsoid (see, e.g., \cite{Blaschke,Santalo,Saint, MP,Petty, MW-3,MR}). Here $V(\cdot)$
denotes the volume (Lebesgue measure) in  $\mathbb{R}^n$ and  the
volume of the Euclidean unit ball $B^n$ in $\Rn$ is solely written
by $\omega_n=\pi^{n/2}/\Gamma(1+n/2)$. A convex body $K$ in $\Rn$ is
a convex compact subset of $\Rn$ with nonempty interior. If the convex
body $K$
contains the origin in its interior, then the polar body $\polar$ of
$K$ is defined  by
\begin{equation*}\polar=\Big\{x\in\mathbb R^n: x\cdot y\le 1\quad\textrm{for
        all}~y\in K\Big\},\end{equation*}  where $x\cdot y$  is the inner
product of $x, y\in \Rn$.  In 1990's, a more general
version of the Blaschke-Santal\'{o} inequality was established by
Lutwak and Zhang \cite{LZ}. It states that, for $K\subset\Rn$ being
a star body and $1\leq p\leq \infty$, one has
\begin{equation}\label{pBS}V(K)V(\Gamma_{p}^{\circ}K)\leq
\omega_n^2,
\end{equation}
with equality if and only if $K$ is an origin-symmetric ellipsoid,
where $\Gamma_{p}^{\circ}K$ denotes the polar body of
the~$L_p$~centroid body $\Gamma_{p}K$ whose support function at
$x\in \Rn$ is given by
\begin{equation}\label{centroid}h_{\Gamma_pK}(x)^p=\frac{1}{c_{n,p}V(K)}\int_K|x\cdot y|^p\,dy,\end{equation} where $\,dy$ is the Lebesgue measure and 
\begin{equation}\label{cnp}c_{n,p}=\frac{\omega_{n+p}}{\omega_2\omega_n\omega_{p-1}},\end{equation} with $\omega_p=\pi^{n/2}/\Gamma(1+p/2)$ for $p>0$. The normalization above is chosen so that $\Gamma_pB^n=B^n$. The
body $\Gamma_{\infty}K$ is to be interpreted as the limit of
$\Gamma_pK$ for $p\rightarrow \infty$. In the case that $K$ is an origin-symmetric convex body and $p=\infty$, $\Gamma_{\infty}^{\circ}K$ coincides with $\polar$, and hence inequality \eqref{pBS}  becomes the Blaschke-Santal\'{o}
inequality \eqref{BS-polar-1}. So inequality
\eqref{pBS} is usually called the $L_p$ Blaschke-Santal\'{o}
inequality.

It can be observed that the definition of the~$L_p$~centroid body is
based on the cosine function (i.e., the inner product). More
specifically, integrating in polar coordinates, for $ x\in \Rn$,
\begin{equation*}\label{e-1}
    h_{\Gamma_pK}(x)^p=\!\frac{n\omega_n}{(n+p)c_{n,p}V(K)}\!\int_{\sphere}|x\cdot
u|^p
\rho_K(u)^{n+p}\,du=\frac{n\omega_n}{(n+p)c_{n,p}V(K)}\big(\mathcal{C}_p\,\rho_K(\cdot)^{n+p}\big)(x),
\end{equation*}
where $\,du$ is the rotation invariant probability measure on the
unit sphere $\sphere$ in $\mathbb{R}^n$ and $\mathcal{C}_p\,\mu$
denotes the $L_p$-cosine transform of a Borel measure $\mu$ defined
on $\sphere$. That is, for $p>0$,
\begin{equation}
   (\mathcal{C}_p\,\mu)(x)=\int_{S^{n-1}}|x\cdot
   u|^p\,d\mu(u), \ \ \ x\in\mathbb{R}^n.\label{p-cosine-transform-2-1}
\end{equation} The $L_p$-cosine transform  provides a very useful
analytical operator for convex geometry and plays a dominating role in applications, see, e.g., the books \cite{Gardner, Koldobsky-FA0, Schneider, HW-1} and references
\cite{BL88,Haberl,HS09,Lon,Lud,Lutwak-5,LYZ04-2,LYZ051,LZ,RZ,Lud02,YY}, among others.

Similar to \eqref{p-cosine-transform-2-1}, $\mathcal{S}_p\,\mu$, the
$L_p$-sine transform of a Borel measure $\mu$ defined on $\sphere$, can be
defined for $p>0$ as follows:
\begin{equation}
   (\mathcal{S}_p\,\mu)(x)=\int_{S^{n-1}} |\textrm{P}_{u^{\perp}}x|^p\,d\mu(u)=\int_{S^{n-1}} [x,
   u]^p\,d\mu(u), \ \ \ x\in\mathbb{R}^n, \label{p-sine-transform-2-1}
\end{equation} where $\textrm{P}_{u^{\perp}}x$ is the orthogonal projection of $x$ onto the $(n-1)$-dimensional subspace  $u^{\perp}$
 perpendicular to $u$, and $[x, u]=\sqrt{|x|^2|u|^2-|x\cdot u|^2}$ with $|x|$ being the Euclidean norm of $x\in \Rn$. Note that $[x, u]$
 represents the sine function due to the Pythagoras theorem, and geometrically $[x, u]$ is the  $2$-dimensional volume of the parallelepiped
 spanned by $x$ and $u$. The  sine, cosine, and Radon  transforms are closely related and have been extensively studied, see, e.g., \cite{Grinberg, Rubin02, Rubin13, Rubin02-2, Ma, Zhang, GW1,GW2,H-S-2019}.  The $L_p$-sine transform (i.e., $m=n-1$) and the $L_p$-cosine transform (i.e., $m=1$) on $S^{n-1}$ are the special cases of the transform $\mathcal{C}_{m,p}\,\mu$ on Grassmann manifolds $\mathbf{G}_{n, m}$ (the set of $m$-dimensional linear subspaces in $\mathbb{R}^n$) \cite{LXZ}. Here, for a Borel measure $\mu$ on $\mathbf{G}_{n, m}$ and $p>0$,
 \begin{equation*}
    (\mathcal{C}_{m,p}\,\mu)(x)=\int_{\mathbf{G}_{n, m}}|\textrm{P}_\xi x|^pd\mu(\xi),\ \ \ x\in\mathbb{R}^n.
 \end{equation*} It can be checked that, when $p\geq 1$, $(\mathcal{C}_{m,p}\,\mu)^{1/p}$ defines a 
 norm on $\mathbb{R}^n$ and hence its unit ball is an origin-symmetric convex body in $\mathbb{R}^n$. Some isoperimetric and reverse isoperimetric inequalities for the unit balls induced by the norm $(\mathcal{C}_{m,p}\,\mu)^{1/p}$ were established in \cite{LXZ}.
Based on the  $L_2$-sine transform,  the first three authors of the present paper proposed  two new sine ellipsoids, and obtained sharp volume inequalities and valuation properties for these sine ellipsoids in \cite{LHX}. Note that these sine ellipsoids are closely related in the Pythagorean relation and duality to  the Legendre ellipsoid and its dual ellipsoid (known as the LYZ ellipsoid) \cite{Lutwak00}.   It is worth mentioning that, when $n=2$, the $L_p$-sine transform coincides with the $L_p$-cosine transform up to a rotation of $\pi/2$. Thus, properties for the $L_p$-cosine transform shall hold for the $L_p$-sine transform on $\mathbb{R}^2$. Unfortunately, this is no long true when $n\geq3$, and the situation is totally different.   The main difficulty for the $L_p$-sine transform is  the lack of affine nature, which makes the analysis on the $L_p$-sine  transform often more challenging.

The main goal of the present paper is to  establish the $L_p$-sine
Blaschke-Santal\'{o} inequality \eqref{main2}, which can be viewed
as the ``sine cousin" of the $L_p$ Blaschke-Santal\'{o} inequality
\eqref{pBS}. Furthermore, a sine version  of  the
Blaschke-Santal\'{o} inequality \eqref{SB} with characterization of
equality is also obtained. Thus, our results are complementary to
the classical studies for the cosine transform.

In Section \ref{sec-3-1}, we introduce a new convex body in
$\mathbb{R}^n$  (i.e., the $L_p$-sine centroid body $\Lambda_pK$ for $K$
being an $L_{n+p}$-star) in terms of the $L_p$-sine transform
\eqref{p-sine-transform-2-1}, whose support function at $x\in \Rn$
is defined by, for $p\geq1$,
\begin{equation}\label{def-p-sine-1}h_{\Lambda_pK}(x)^p=\frac{n\omega_n}{(n+p)\widetilde{c}_{n,p}V(K)}\big(\mathcal{S}_p\,\rho_K(\cdot)^{n+p}\big)(x)=\frac{1}{\widetilde{c}_{n,p}V(K)}\int_K[x,y]^p\,dy,\end{equation}
where
\begin{equation}\label{tcnp}\widetilde{c}_{n,p}=\frac{(n-1)\omega_{n-1}\omega_{n+p-2}}{(n+p)\omega_n\omega_{n+p-3}}.\end{equation}
In particular, $\Lambda_pB^n=B^n$ (see
\eqref{B}). Note that $\Lambda_pK$ is well-defined since
$h_{\Lambda_pK}$ is a convex function (see
 Section \ref{sec-3-1} for details).
 In the special case when $K$ is a star body in $\Rn$,  $\Lambda_pK$ was first proposed in
\cite{LHX}.
Motivated by the limit of $\Lambda_p^{\circ} K$  as $p\rightarrow
\infty$,  the sine polar body of $K$, denoted by $\ksp$, is
introduced. Namely, for a subset $K\subset \Rn$, we let
$$\ksp=\Big\{x\in\Rn: [x, y]\le 1\ \ \mathrm{for \ all}\ y\in
K\Big\}.$$ The sine polar body $\ksp$ differs from $\polar$
mainly with $[x, y]$ replacing $x\cdot y$. Hence many properties for
$\ksp$ are similar to those for $\polar$, for instance,
$L^{\diamond}\subseteq \ksp$ if $K\subseteq L$,
$(\ball)^{\diamond}=\ball$, and $(cK)^{\diamond}=c^{-1}\cdot \ksp$
 for any $c>0$. However, some properties for $\ksp$ are completely different from those for $\polar$.
For instance, the volume product $V(K)V(\polar)$
is~$\textrm{GL}(n)$-invariant (i.e., $V(\phi K)V((\phi
K)^{\circ})=V(K)V(\polar)$ for any invertible linear map $\phi$ on
$\Rn$), but the volume product $V(K)V(\ksp)$ is in general not
$\textrm{GL}(n)$-invariant (except for $n=2$). The
non-$\textrm{GL}(n)$-invariance for $V(K)V(\ksp)$ does bring extra
difficulty in  characterizing the equality for the sine
Blaschke-Santal\'{o} inequality \eqref{SB} in Section
\ref{BS-Sine-2-2-1}.  On the other hand, the sine bipolar property
$K=K^{\diamond\diamond}$, where
$K^{\diamond\diamond}=(\ksp)^{\diamond}$, does not hold in general
for all convex bodies. Indeed, the definition of $\ksp$ indicates
that $\ksp$ is  formed by the intersection of closed solid cylinders
(see \eqref{cy}). This brings our attention to a special class of
convex bodies $\Cyk$ consisting of all origin-symmetric convex
bodies generated by the intersection of origin-symmetric closed
solid cylinders. The sine bipolar property $K=K^{\diamond\diamond}$ holds for all $K\in \Cyk$, see Proposition \ref{prop-p-s} (v). Proposition \ref{prop-sine-polar-1} (ii) also shows
that $K^{\diamond\diamond}$ is just the cylindrical hull of $K$ --
the intersection of all origin-symmetric closed solid cylinders
containing $K$. The cylindrical hull of $K$ is a notion analogous to
the convex hull, and any convex body $K\in \Cyk$ must be equal to
the cylindrical hull of itself (see Proposition
\ref{prop-sine-polar-1} (i)). Moreover, in Section \ref{sec-3-1}, we  develop
many notions in $\mathcal{C}_{e}^n$, such as the cylindrical support
function, the supporting cylinder, and the cylindrical Gauss image.
These notions are completely parallel to their classical
counterparts, e.g., the support function, the supporting hyperplane,
and the Gauss image.

\begin{figure}[htbp]
\centering
\begin{minipage}[t]{0.48\textwidth}
\centering
\includegraphics[width=4cm]{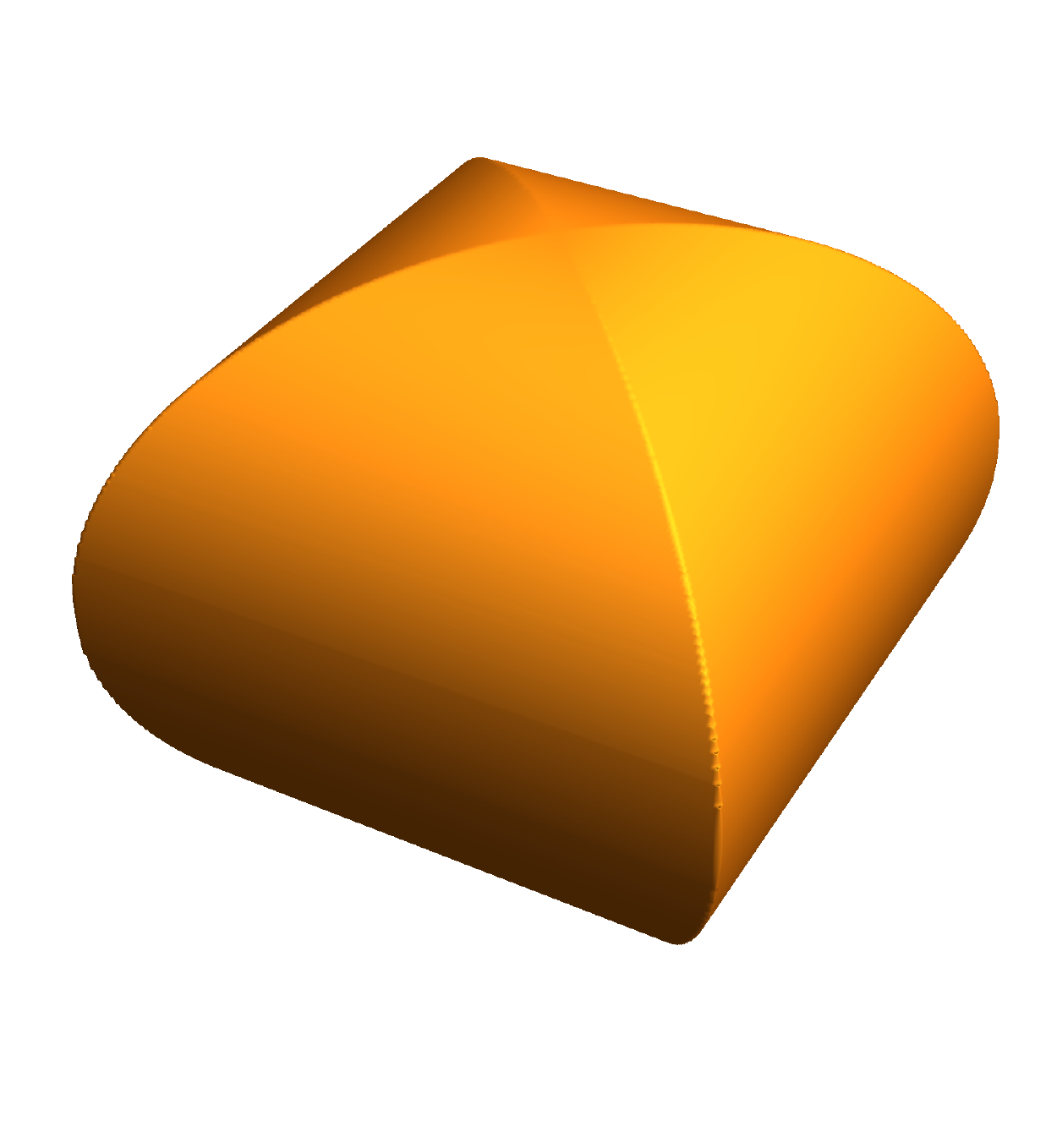}
\end{minipage}
\begin{minipage}[t]{0.48\textwidth}
\centering
\includegraphics[width=4cm]{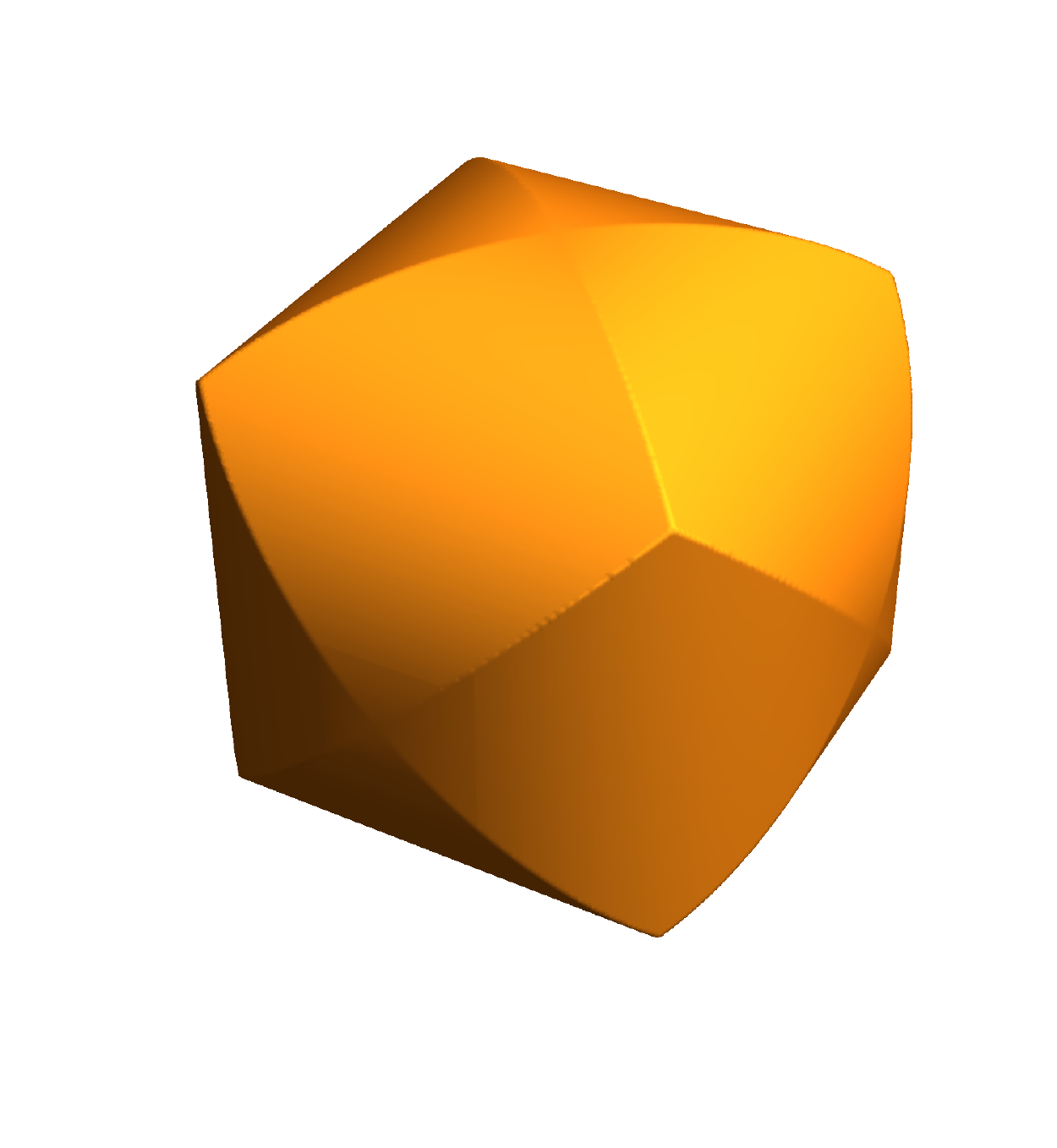}
\end{minipage}
\caption{The left one is the bicylinder (``Mou He Fang Gai" in
Chinese) and the right one is a tricylinder in $\mathbb R^3$.  More
figures can be found in e.g., \cite{Paul-2003, H-I}.} \label{fig-1}
\end{figure}

 Last but not the least, we would like to mention that the notion $\Cyk$ is quite important. A basic example of $\Cyk$ is the intersection of two closed solid cylinders with equal radius at right
 angles. The volume of this bicylinder (also known as ``Mou He Fang Gai" in China) in $\mathbb R^3$ was first studied by Archimedes more than 2200 years ago and independently by the old Chinese mathematicians  Hui Liu,  Chongzhi Zu and Geng  Zu  more than 1500 years ago (see, e.g., \cite{Jan-2002, TK-1972}). An amazing application of the bicylinder is the calculation of the volume of $3$-dimensional
balls by these old Chinese mathematicians \cite{TK-1972}. The solids
obtained as the intersection of two or three cylinders of equal
radius at right angles are better known as the Steinmetz solids in
Europe, see Figure \ref{fig-1}. Applications of the Steinmetz solid
are common in real life, such as the T-adapter and the circular
pipes joining.  The intersection of cylinders also naturally appears
in the study of crystal (see, e.g., \cite{I-O-A, Moore-1974}).

In Section \ref{sec-3-2}, we will prove the following $L_p$-sine
Blaschke-Santal\'{o} inequality: {\em let  $p\geq 1$ and $K\subset
\Rn$ be a star body. Then
\begin{equation}\label{main2}V(K)V(\plp)\leq \omega_n^2,\end{equation}
with equality if and only if $K$ is an origin-symmetric ellipsoid
when $n=2$, and is an origin-symmetric ball when $n\geq3$.} For
$p=2$, this inequality was also established by the first three
authors \cite[Theorem 4.1]{LHX}. By letting $p\rightarrow \infty$,
the $L_p$-sine Blaschke-Santal\'{o} inequality will lead to the
following sine Blaschke-Santal\'{o} inequality: {\em let $K\subset
\Rn$ be a star body. Then
 \begin{equation}\label{SB}
V(K)V(\ksp) \leq  \omega_n^{2},
\end{equation} with equality if and only if $K$ is an origin-symmetric ellipsoid
when $n=2$, and is an origin-symmetric ball when $n\geq3$.}

Therefore, inequalities \eqref{main2} and \eqref{SB} can be viewed
as the ``sine cousin" of inequalities \eqref{pBS} and
\eqref{BS-polar-1}. However, the approximation as
$p\rightarrow \infty$ does not yield the equality conditions of
inequality \eqref{SB}. To characterize these equality conditions
requires a lot of work mainly because of the extra difficulty
arising from the non-$\textrm{GL}(n)$-invariance of $V(K)V(\ksp)$
for $n\geq 3$. In this sense, equality conditions in inequality
\eqref{BS-polar-1} and its sine counterpart \eqref{SB} are different
for $n\geq 3$. However, inequalities \eqref{BS-polar-1} and \eqref{SB}
are in fact equivalent if $K$ is an origin-symmetric convex body in
$\mathbb R^2$. Indeed, it will be shown in Section \ref{sec-3-2}
that under this circumstance $\ksp$ is just a rotation of $\polar$
by an angle of $\pi/2$. A similar situation also happens to
 \eqref{pBS} and \eqref{main2}.

We also  mention that information theory
 (or probability theory) and convex geometry are closely related (see, e.g., \cite[Section 14]{Gardner02}).
In recent years there was a growing body of works in this direction,
see, e.g., \cite{BM,ALYZ-14,FM, GLYZ-02,  LYZ-duke, LYZ04, LYZ05,
LYZ07,LYZ12,LYZ13, Ng, MMX}. For example, an important application
of the $L_p$ Blaschke-Santal\'{o} inequality \eqref{pBS} to
information theory leads to the celebrated $L_p$ moment-entropy
inequality (see \eqref{mm}) due to Lutwak, Yang, and Zhang \cite{LYZ04}.
Similarly, the ``sine cousin" of $L_p$ moment-entropy inequality (see
\eqref{inf}) can be deduced by using the $L_p$-sine
Blaschke-Santal\'{o} inequality. We believe that the applications of
the $L_p$-sine transform to information theory (or probability
theory) will likely produce more interesting and useful results for
applications in many areas.

\section{Preliminaries}\label{section-b-2}
Throughout the paper, $\Rn$ denotes $n$-dimensional Euclidean space
($n\geq 2$) and the term subspace means a linear subspace. As
usual,~$x\cdot y$~denotes the inner product of $x,y\in\Rn$
and~$|x|$~the Euclidean norm of~$x$. The unit sphere $\{x\in
\Rn:|x|=1\}$ is denoted by $\sphere$ and the unit ball $\{x\in
\Rn:|x|\leq 1\}$ by $\ball$. Furthermore, we denote by $o$ the origin
in $\Rn$ and $\bar{x}=x/|x|\in \sphere$ for all $x\in \Rn\setminus
\{o\}$. Note that the cosine of the angle between nonzero vectors
$x$ and $y$ is indeed given by $\bar{x}\cdot \bar{y}$. Thus, for $x,
y\in \Rn\setminus \{o\}$, \begin{equation} \label{cos}|x\cdot y|
=|x|\cdot |y|\cdot |\bar{x}\cdot \bar{y}| =|x|\cdot
|\textrm{P}_{x}y|,\end{equation}  where $\textrm{P}_{x}y$ is the
orthogonal projection of $y$ onto the $1$-dimensional subspace
containing $x$. On the other hand, the sine of the angle between
nonzero vectors $x$ and $y$ can be easily calculated through the
formula $\sqrt{1-|\bar{x}\cdot \bar{y}|^2}$. Corresponding to
\eqref{cos}, for $x, y\in \Rn\setminus \{o\}$, let
\begin{equation} \label{cal-sin} [x,y] :=|x|\cdot
|y|\cdot\sqrt{1-|\bar{x}\cdot \bar{y}|^2}=\sqrt{|x|^2|y|^2-|x\cdot
y|^2} =|x|\cdot\sqrt{|y|^2-|\textrm{P}_{x}y|^2}=|x|\cdot
|\textrm{P}_{x^{\perp}}y|,
\end{equation} where $\textrm{P}_{x^{\perp}}y$ is the orthogonal projection of $y$ onto the $(n-1)$-dimensional subspace $x^{\perp}$ perpendicular to $x$.
In other words, $[x, y]$ represents the  $2$-dimensional volume of
the parallelepiped spanned by $x$ and $y$.  Of course, one can switch the
roles of $x, y$ to get $[x,y]=[y,x]=|y|\cdot
|\textrm{P}_{y^{\perp}}x|$. Moreover, if one or both of $x$ and $y$
are zero vectors, one can simply let $[x, y]=0$.

A convex body in $\Rn$ is a convex compact subset of $\Rn$ with
nonempty interior. By $\cKon$ we mean the class of all convex bodies
with the origin $o$ in their interiors. For each $K\in \cKon$, one
can define its support function $h_K: \Rn \rightarrow (0,\infty)$ by
\begin{equation}\label{support}
h_K(x)=\max \{x\cdot y: y\in K\},\ \ \ x\in \Rn.
\end{equation} A convex body $K\in \cKon$ is said to be origin-symmetric if $-x\in K$ for all $x\in K$. Denote by $\cKe$ the subclass of all origin-symmetric elements of
$\cKon$.

Denote by $\textrm{GL}(n)$ the group of all invertible linear
transforms on $\Rn$ and $\textrm{SL}(n)$ the subgroup of
$\textrm{GL}(n)$ whose determinant is $1$. By $\textrm{O}(n)$ we
mean  the orthogonal group of~$\Rn$, that is, the group of linear
transformations preserving the inner product. For each $\phi\in
\textrm{GL}(n)$, let $\phi^t$ and $\phi^{-1}$ be the transpose and
inverse of $\phi$, respectively. For $\phi\in \textrm{GL}(n)$ and
$x\in \Rn$, let $\phi K=\{\phi y: y\in K\}$ and then
\begin{equation}\label{affh}h_{\phi K}(x)=h_{K}(\phi^t x).\end{equation}

 Throughout the paper, $\,du$ denotes the rotation invariant probability measure on
$\sphere$.  Let $q> 0$. We say that a function $\rho: \sphere
\rightarrow [0, \infty)$ is $q$-integrable on $\sphere$  if $\rho$
is measurable and $$\int_{\sphere} \rho(u)^q\,du<\infty.$$ Denote by
$L^{q}(S^{n-1})$ the set of all $q$-integrable functions on
$\sphere$. As $\int_{\sphere}\,du=1$, one can easily obtain, by
using the H\"{o}lder inequality, that if $0<q_1<q$, then
\begin{equation}\label{compare-space} L^q(\sphere)\subseteq
L^{q_1}(\sphere).\end{equation}

A set $L\subset\Rn$ is said to be star-shaped about the origin if
$o\in L$ and
$$L=\big\{ru: 0\le r \le \rho(u)\ \ \textrm{for} \ \ u\in S^{n-1}\big\},$$
where $\rho: \sphere\rightarrow [0, \infty)$  is a nonnegative
function. Such a function $\rho$ is called the radial function of
$L$, which will be more often written as $\rho_L.$ The radial
function can be extended to $\rho_L: \Rn\setminus\{o\} \rightarrow
[0, \infty)$ by $\rho_{L}(ru)=r^{-1}\rho_L(u)$ for $r>0$ and $u\in
\sphere$, which gives
\begin{equation}\label{radial}\rho_L(x)=\max\big\{\lambda\ge 0:\lambda
x\in L\big\},\ \ \ x\in \Rn\setminus\{o\}.\end{equation} Let $p>0$. A star-shaped set $L$ in
  $\Rn$ is said to be an \textit{$L_{n+p}$-star}  generated by
$\rho\in L^{n+p}(S^{n-1})$ if $\rho_L=\rho$. Clearly, it follows
from \eqref{radial} that for $c>0$ and $x\in \mathbb
R^n\backslash\{o\}$,
\begin{equation}\label{rad2}
    \rho_{cL}(x)=c\rho_L(x).
\end{equation}
For $\phi\in \textrm{GL}(n)$ and $x\in \mathbb R^n\backslash\{o\}$,
one has $\rho_{\phi L}(x)=\rho_{L}(\phi^{-1}x)$.

A set $L\subset \Rn$ is  called a star body if $L$ is star-shaped
and $\rho_L$ restricted on $\sphere$ is positive and continuous.
Denote by $\cS^n$ the class of all star bodies.  Clearly,
$\cKon\subset \cS^n$, and $L$ is an $L_{n+p}$-star  for all $p>0$ if
$L\in \cS^n$. Moreover, the Minkowski
functional~$\|\cdot\|_L:\Rn\rightarrow [0,\infty)$~of $L\in \cS^n$
is defined by
\begin{equation}\label{Minkow}\|x\|_L=\min\{\lambda\ge 0:
x\in \lambda L\},\ \ \ x\in \Rn.\end{equation}

Recall that the polar body $\polar$ of $K\in \cKon$ is defined by
\begin{equation}\label{pol}\polar=\Big\{x\in\mathbb R^n: x\cdot y\le 1\quad\textrm{for
all}~y\in K\Big\}.\end{equation} It can be easily verified that for
$K\in \cKon$,\begin{eqnarray}
    (c K)^{\circ}\!\!\!&=&\!\!\!c^{-1}\polar, \ \ \ \mathrm{for \ any}\ c>0, \label{cm}
\\ (\phi K)^{\circ}\!\!\!&=&\!\!\!\phi^{-t}\polar,\ \ \ \mathrm{for \ any}\ \phi\in
\textrm{GL}(n). \label{19}
\end{eqnarray} Moreover, it follows from \eqref{support}, \eqref{radial}, \eqref{Minkow}, and \eqref{pol} that if $K\in \cKon$,
\begin{equation}\label{rh}
   \|x\|_{\polar}=1/\rho_{\polar}(x)=h_{K}(x),\ \ \ \mathrm{for\ any}\ x\in\Rn\setminus\{o\}.
\end{equation} In particular, $(\polar)^\circ=K$ for any $K\in \cKon$. Note that if $K\in \cKe$ is an origin-symmetric convex body
in~$\Rn$,  \eqref{pol} can also be rewritten as
\begin{equation}\label{polar}
    \polar=\{x\in\mathbb R^n: |x\cdot y|\le 1\quad\textrm{for all}~y\in
    K\}.
\end{equation}

 Let $p>0$. For $L$ being an $L_{n+p}$-star, its volume is given
 by
 \begin{equation}\label{volume-12-5}
V(L)=\omega_n\int_{\sphere}\rho_L(u)^n\,du.
\end{equation} By \eqref{compare-space}, for an $L_{n+p}$-star $L$,
it has a finite volume as $\rho_L\in L^{n+p}(\sphere).$  Moreover,
the $L_{-p}$ dual mixed volume \cite{LYZ04} of an $L_{n+p}$-star $L$ and a
star body $L'\in \cS^n$ can be formulated by
\begin{equation}\label{dual}
\widetilde{V}_{-p}(L,
L')=\omega_n\int_{\sphere}\rho_L(u)^{n+p}\rho_{L'}(u)^{-p}\,du.
\end{equation} In particular, for $L\in \cS^n$ and $p>0$,  one has \begin{equation}\label{e-21}
   \widetilde{V}_{-p}(L, L)=V(L).
\end{equation} Following from the H\"{o}lder inequality, one can easily obtain the following dual $L_{-p}$ Minkowski inequality \cite{LYZ04}:  for
$p>0$,
\begin{equation}\label{dm}
     \widetilde{V}_{-p}(L, L')^n \geq V(L)^{n+p}V(L')^{-p}
\end{equation} holds for any $L_{n+p}$-star  $L$  and any
$L'\in \cS^n$,  with equality if and only if there exists a constant
$c>0$ such that  $\rho_{L'}(u)=c\cdot\rho_L(u)$ for almost all $u\in
\sphere$ with respect to the spherical measure of $\sphere$.

\section{The $L_p$-sine centroid body, the sine polar body, and the cylindrical hull}\label{sec-3-1}
Let $p\ge 1$ and $n\ge 2$. We first show that, under the assumption
that $K\subset\Rn$ is an $L_{n+p}$-star with $V(K)>0$, the set
$\Lambda_pK$ defined in \eqref{def-p-sine-1} is well-defined.
Indeed, we can get  $ h_{\Lambda_pK} \in (0, \infty)$ on
$\Rn\setminus\{o\}$. To see this, for any $x\in \Rn\setminus\{o\}$
and any $y\in K$,   (\ref{cal-sin}) implies $[x, y]\leq |x|\cdot
|y|$ and hence
    \begin{align*}
    h_{\Lambda_pK}(x)^p &=\frac{1}{\widetilde{c}_{n,p}V(K)}\int_K[x,y]^p\,dy\\ &\leq \frac{1}{\widetilde{c}_{n,p}V(K)}\int_K|x|^p|y|^p \,dy\\& =\frac{n\omega_n |x|^p}{(n+p)\widetilde{c}_{n,p}V(K)} \int_{\sphere}\rho_K(u)^{n+p}\,du<\infty.  \end{align*}
    On the other hand, as $V(K)>0$, then $K$ cannot be concentrated on a $1$-dimensional space, and  $[x, y]>0$  for almost all $y\in K$
    (with respect the Lebesgue measure). Hence, $h_{\Lambda_pK}(x)>0$ for all $x\in \Rn\setminus\{o\}$.

\vskip 2mm   It can be easily verified that  $h_{\Lambda_pK} (tx)=t\cdot h_{\Lambda_pK}(x)$ for all $t\geq 0$ and $x\in \Rn$. Moreover, $h_{\Lambda_pK}(x)$ is a convex function, namely for all $\tau\in [0, 1]$ and $x, x'\in \Rn$, one has
$$h_{\Lambda_pK}(\tau x+(1-\tau)x')\leq \tau\cdot h_{\Lambda_pK}(x)+(1-\tau)\cdot h_{\Lambda_pK}(x').$$
To see this, for $y\in \Rn\setminus\{o\}$, one has $[x, y]=|y|\cdot
|\textrm{P}_{y^{\perp}}x|$ and hence
\begin{eqnarray} [\tau x+(1-\tau)x', y]\!\!&=&\!\!|y|\cdot \big|\textrm{P}_{y^{\perp}}\big(\tau x+(1-\tau)x'\big) \big|
\nonumber  \\\!\!&=&\!\!|y|\cdot \big|\tau  \textrm{P}_{y^{\perp}}x
+(1-\tau)\textrm{P}_{y^{\perp}}x' \big| \nonumber
\\\!\!&\leq&\!\! \tau|y|\cdot \big|\textrm{P}_{y^{\perp}}x\big|
+(1-\tau)|y|\cdot \big|\textrm{P}_{y^{\perp}}x' \big| \nonumber
\\\!\!&=&\!\!\tau [x, y]+(1-\tau)[x', y]. \label{conv--1}
\end{eqnarray} Clearly, \eqref{conv--1} is also valid for~$y=o$. It
then follows from the Minkowski inequality that
\begin{eqnarray*}h_{\Lambda_pK}(\tau
x+(1-\tau)x')\!\!&=&\!\!\bigg(\frac{1}{\widetilde{c}_{n,p}V(K)}\int_K[\tau
x+(1-\tau)x',y]^p\,dy\bigg)^{1/p} \nonumber\\ \!\!&\leq&\!\!
\bigg(\frac{1}{\widetilde{c}_{n,p}V(K)}\int_K\Big(\tau [x,
y]+(1-\tau)[x',y]\Big)^p\,dy\bigg)^{1/p}   \nonumber\\
\!\!&\leq&\!\! \tau\cdot
\bigg(\frac{1}{\widetilde{c}_{n,p}V(K)}\int_K [x,
y]^p\,dy\bigg)^{1/p} +(1-\tau)\cdot
\bigg(\frac{1}{\widetilde{c}_{n,p}V(K)}\int_K [x', y]^p\,dy\bigg)^{1/p} \nonumber\\
\!\!&=&\!\!\tau\cdot h_{\Lambda_pK}(x)+(1-\tau)\cdot
h_{\Lambda_pK}(x').
\end{eqnarray*} Furthermore, $\Lambda_pK$ is clearly origin-symmetric. We conclude that $\Lambda_pK\in \cKe$ for $p\geq 1$.

Denote by $\plp$ the polar body of $\Lambda_{p}K$.  It follows from  \eqref{rh},
\eqref{def-p-sine-1}, and the  polar coordinate that
\begin{eqnarray}\label{Po}
    \rho_{\plp}(x)^{-p}=\frac{n\omega_n}{(n+p)\widetilde{c}_{n,p}V(K)}\int_{\sphere}[x, u]^p\rho_K(u)^{n+p}\,du,\ \ \ x\in
    \Rn\setminus\{o\}.
\end{eqnarray}
By \eqref{rad2} and \eqref{Po}, one has
\begin{equation}\label{c}
    \Lambda_p^{\circ}(cK)=c^{-1}\plp,\ \ \ \mathrm{for \ any}\ c>0.
\end{equation}

If $K\in \cS^n$ is a star body in $\mathbb{R}^n$, then $K$ is an $L_{n+p}$-star
for all $p\geq 1$ with $V(K)>0$. Thus, by
\eqref{def-p-sine-1} and Stirling's formula,
\begin{equation*}
    h_{\Lambda_{\infty}K}(x)=\lim_{p\rightarrow\infty}h_{\Lambda_{p}K}(x)=\max_{y\in
   K}\,[x, y],\ \ \  \mathrm{for\ any}\ x\in \Rn.
\end{equation*} Again it is easy to check that $\Lambda_{\infty}K\in \cKe$. Denote by $\pli$ the polar body of $\Lambda_{\infty}K$.  It follows from \eqref{rh} that \begin{eqnarray}
 \label{Po-1}
    \|x\|_{\pli} =h_{\Lambda_{\infty}K}(x)=\max_{y\in
   K}\,[x, y],\ \ \ \mathrm{for\ any}\ x\in \Rn.
\end{eqnarray}
Formula \eqref{Minkow} for the Minkowski functional  yields that,
for any star body $L\in \cS^n$, one must have $L=\big\{x\in \Rn: \
\|x\|_{L}\leq 1\big\}.$ Together with \eqref{Po-1}, one gets, if
$K\in \mathcal{S}^n$,
\begin{equation*}
   \pli=\Big\{x\in \Rn: \ \|x\|_{\pli}\leq 1\Big\}=\Big\{x\in \Rn: \ \max_{y\in
   K}\,[x, y]\leq 1\Big\}.
\end{equation*}

This motivates our definition for the
sine polar body. \bd \label{s-p-b} Let $K$ be a subset in $\Rn$ with
$n\geq 2$. Define $\ksp$,  the {\em sine polar body} of $K$, to be
\begin{equation}\label{sc}\ksp=\Big\{x\in\Rn: [x, y]\le 1\ \ \mathrm{for
\ all}\ y\in
    K\Big\}.\end{equation}\ed
Clearly,  for $K\in \cS^n$, one has
\begin{equation}\label{dinf}\ksp=\pli. \end{equation}  If $K\in \cKe$,
one sees that the major difference of \eqref{polar} and \eqref{sc}
is in fact the replacement of $|x\cdot y|$ by $[x, y]$, and hence
the newly defined sine polar body $\ksp$ can be viewed as the sine
counterpart of $\polar$.

Throughout the paper, a closed solid cylinder $C^{-}(u,
r)\subset\Rn$ with axis being the line $\{tu: t\in \R\}$ for $u\in
\sphere$ and base radius $r>0$ is a subset of $\Rn$ of the following
form:
\begin{equation*}\label{solid-cyl} C^{-}(u, r)=\big\{x\in \Rn: [x, u]\leq
r\big\}=\big\{x\in \Rn: |\textrm{P}_{u^{\perp}}x|\leq
r\big\}.\end{equation*} Due to \eqref{polar}, one sees that if $K\in
\cKe$, then $\polar$ is obtained by the intersection of slabs.
However, it follows from \eqref{cal-sin} and \eqref{sc} that
\begin{eqnarray} \ksp\!\!&=&\!\!\Big\{x\in\Rn: |y|\cdot |\textrm{P}_{y^{\perp}}x|\le 1\ \ \mathrm{for
\ all}\ y\in
    K\setminus\{o\}\Big\}
\nonumber  \\\!\!&=&\!\!\Big\{x\in\Rn: |\textrm{P}_{y^{\perp}}x|\le
\frac{1}{|y|}\ \ \mathrm{for \ all}\ y\in
    K\setminus\{o\}\Big\}.\label{cy}
\end{eqnarray}
This means that $\ksp$ is indeed formed by the intersection of
closed solid cylinders (with axis being the line $\{ty: t\in \R\}$
and base radius being $|y|^{-1}$ for $y\in
    K\setminus\{o\}$). As the polar body plays
fundamental roles in convex geometry, such as in the celebrated
Blaschke-Santal\'{o} inequality \eqref{BS-polar-1} and many other
affine isoperimetric inequalities (see, e.g., \cite{Lut}), we expect
that the sine polar body will play roles similar to its ``cosine
cousin" in applications. In particular, as a start, the sine
Blaschke-Santal\'{o} inequality for the sine polar body will be
established in Section \ref{BS-Sine-2-2-1}.

In view of Definition \ref{s-p-b} and the nature of the sine polar
body, it is worth to investigate $\Cyk\subseteq\cKe$,  the family of
origin-symmetric convex bodies generated by the intersection of
origin-symmetric closed  solid cylinders.  Typical examples of
$\Cyk$ are the Steinmetz solids, see Figure \ref{fig-1}. As pointed out by one of the referees, our sine polarity turns out to be a special case of the notion of polarity with respect to a
general function $c(\cdot, \cdot)$  introduced by Artstein-Avidan,
Sadovsky, and Wyczesany \cite{ASW} (namely, by letting $c(x, y)=-[x, y]$ and $t=-1$).  Consequently, $\Cyk$ is a type of the $c$-class of sets in \cite{ASW}. These concepts were further studied in \cite{ASW-1}.

It is indeed a
surprise that to generate a convex body $K\in \Cyk$, as few as only
two closed solid cylinders whose axes are not parallel to each other
would be enough. In other words, one needs the directions of the
axes of the closed solid cylinders generating $K\in \Cyk$ being not
concentrated on an $1$-dimensional subspace. This result is
summarized in the following proposition.

 \bp \label{non-conc-cyl}  Let  $r_1, r_2>0$ be  two  constants  and directions $u, v\in \sphere$ such that $u$ is not parallel to $v$. Then $C^{-}(u, r_1)\cap C^{-}(v, r_2)\in \Cyk.$  \ep
\begin{proof} Let $K=C^{-}(u, r_1)\cap C^{-}(v, r_2)$. To prove $K\in \Cyk$, it would be enough to prove $K\in \cKe$.

 First of all, it is clear that $K$ is an origin-symmetric closed convex set with nonempty interior.
 Suppose that $K$ is unbounded. Then there exists a direction $w\in\sphere$
 such that $\{\tau w: \tau>0\} \subseteq K$. In particular, $[\tau w, u]\leq r_1$ for all $\tau>0$.
 This implies  $0\leq [w, u]\leq r_1/\tau\rightarrow 0$ as $\tau\rightarrow \infty$. Thus, $[w,u]=|\textrm{P}_{u^{\perp}}w|=0$, and  $u$ is parallel to  $w$. Similarly,  $v$ is parallel to $w$, and thus $u$ is parallel to $v$, a contradiction. So  $K=C^{-}(u, r_1)\cap C^{-}(v, r_2)$ is bounded and then $K\in \cKe$ as desired. \end{proof}

Basic properties for the sine polar body can be listed as follows.

    \bp \label{prop-p-s} Let $K$ be a subset in $\Rn$ with $n\geq 2$. Then the following properties hold.
    \vskip 2mm \noindent \emph{(i)}  Suppose $L\subseteq \Rn$ such that $K\subseteq L$. Then $L^{\diamond}\subseteq \ksp$.
  \vskip 2mm \noindent \emph{(ii)}  The sine polar body of $\ball$ is again $\ball$, that is,  $(\ball)^{\diamond}=\ball$.
 \vskip 2mm \noindent \emph{(iii)} For any $O\in \textrm{O}(n)$ and $c>0$, one has, $(OK)^{\diamond}=O\ksp$ and $(cK)^{\diamond}=c^{-1}\cdot \ksp$.
 \vskip 2mm \noindent \emph{(iv)}  If $K\subset \Rn$ is bounded and is not concentrated on any $1$-dimensional subspace, then   $\ksp\in \Cyk$.
  \vskip 2mm \noindent \emph{(v)} Let $K^{\diamond \diamond}=(\ksp)^\diamond.$ Then $K\subseteq K^{\diamond \diamond}$. Moreover, $K= K^{\diamond \diamond}$ for all $K\in \Cyk$.
   \vskip 2mm \noindent \emph{(vi)} For any $K\in \Cyk$, one has \begin{equation}\label{compare-polar-diamond} V(K^{\diamond})\leq V(K^{\circ}). \end{equation}
  \ep

We omit the proofs of properties (i)-(v), since they are easy to get.  Properties (vi) is vital to prove the sine
Blaschke-Santal\'{o} inequality, which will be proved in the end of
this section. Note that properties (i) and (v) have also appeared in
\cite[Lemma 3.10]{ASW} for $c$-polarity.

    Recall that $K^{\circ\circ}$ is the smallest convex body that contains a compact subset $K\subset \Rn$ with the origin in its interior, namely, $K^{\circ\circ}$ is the convex hull of $K$. Motivated by Proposition \ref{prop-p-s},  the following {\em cylindrical hull} of a bounded set may be proposed.

   \bd\label{cyl-hull} Let $E\subset\Rn$ be bounded and not concentrated on any $1$-dimensional subspace. The cylindrical hull of $E$, denote by $\mathrm{cyl}(E)$, is a convex body in $\Rn$ of the following form: $$\mathrm{cyl}(E)=\bigcap\big\{C: \ C\ \mathrm{is\ a \ closed \ solid\ cylinder\ such\ that } \ E\subset C\big\}.$$
   \ed
If $E\subset\Rn$ is bounded and is not concentrated on any
$1$-dimensional subspace, then $E\subseteq \mathrm{cyl}(E)$ and
$\mathrm{cyl}(E)\in \Cyk.$ Indeed, as $E$ is not concentrated on any
$1$-dimensional subspace, there must be two nonparallel vectors $x,
y\in \Rn\setminus\{o\}$ such that $x,y\in E$. It further yields that
all the base radii of closed solid cylinders containing $E$ are
positive. Therefore, $\textrm{cyl}(E)\in \Cyk$. Since
$\textrm{cyl}(E)$ is a convex body in $\Rn$ containing $E$, one also
has $\mathrm{cyl}(E)\supseteq \mathrm{conv}(E)$, as the convex hull
$\mathrm{conv}(E)$ is the smallest convex body containing $E$.

The following result asserts that $K^{\diamond\diamond}$ is exactly
the cylindrical hull of $K$. Again, the cylindrical hull turns out to be a special
case of the $c$-envelope defined in \cite{ASW}.
\bp\label{prop-sine-polar-1} Let $K\subset\Rn$ be bounded and not
concentrated on any $1$-dimensional subspace. Then the following
results hold. \vskip 2mm \noindent \emph{(i)} $K\in \Cyk$ if and
only if $K=\mathrm{cyl}(K)$. \vskip 2mm \noindent \emph{(ii)}
 $K^{\diamond\diamond}=\mathrm{cyl}(K).$ \ep

\begin{proof} (i) Let $K\in \Cyk$. Note that $K\subseteq \mathrm{cyl}(K)$. Hence, $K=\mathrm{cyl}(K)$ follows immediately if $K\supseteq \mathrm{cyl}(K)$ is verified. This is an easy consequence of $K\in \Cyk$ as $K$ is an intersection of closed solid cylinders (of course) containing $K$.

\vskip 2mm \noindent (ii) First of all, it follows from Proposition
\ref{prop-p-s} (iv) that $K^{\diamond}\in \Cyk$ and
$K^{\diamond\diamond}\in \Cyk$. Together with $K\subseteq
K^{\diamond\diamond}$ due to Proposition \ref{prop-p-s} (v), one
gets that $K^{\diamond\diamond}$ is an intersection of closed solid
cylinders containing $K$. This implies that
$\mathrm{cyl}(K)\subseteq K^{\diamond\diamond}.$

Now let us prove $\mathrm{cyl}(K)\supseteq K^{\diamond\diamond}.$ By
Proposition \ref{prop-p-s} (i), the fact that $K\subseteq
\mathrm{cyl}(K)$ implies $\ksp\supseteq \mathrm{cyl}(K)^{\diamond}$.
From Proposition \ref{prop-p-s} (i), (v) and the fact that
$\mathrm{cyl}(K)\in \Cyk$, we further get
$K^{\diamond\diamond}\subseteq
\mathrm{cyl}(K)^{\diamond\diamond}=\mathrm{cyl}(K)$ as desired.
\end{proof}

Recall that any convex body $K\in \cKe$ can be formed by
$$K=\bigcap_{u\in \sphere} \Big\{x\in \Rn: |x\cdot u|\leq
h_K(u)\Big\}.$$ Note that $h_K(u)=\max_{y\in K} u\cdot y$ is the
distance from the origin to the tangent hyperplane of $K$ at the
direction $u\in \sphere$. Likewise, we can propose a definition for
 the cylindrical support function.

\bd \label{cyl-supp-fun-def} The cylindrical support
function $c_K: \Rn\rightarrow [0, \infty)$ of $K\in \Cyk$  is defined by
\begin{equation}\label{cyl-supp-func} c_K(x)=\max_{y\in K}\,
[x,y],\ \ \ x\in \Rn.\end{equation} \ed The function $c_K$ is an
even function and has positive homogeneity of degree $1$ in the
sense that $c_K(rx)=r\cdot c_K(x)$ for all $r\geq 0$ and $x\in \Rn$.
Moreover, $c_K$ is convex, that is
$$c_K(\tau x+ (1-\tau)x')\leq \tau c_K(x)+(1-\tau)c_K(x')$$ for any
$x, x'\in \Rn$ and $\tau\in [0, 1]$. Indeed, by \eqref{conv--1}, one
has
\begin{eqnarray*} c_K(\tau x+ (1-\tau)x')\!\!&=&\!\!\max_{y\in K}\,[\tau
x+(1-\tau)x', y] \\ \!\!&\leq&\!\! \max_{y\in K}\Big( \tau [x,
y]+(1-\tau)[x', y]\Big)\\ \!\!&\leq&\!\! \tau \max_{y\in K}\,[x,
y]+(1-\tau)\max_{y\in K}\,[x', y]\\ \!\!&=&\!\!\tau
c_K(x)+(1-\tau)c_K(x').
\end{eqnarray*}

Denote by $\partial K$  the boundary of $K$. The  cylindrical
support function can be used to form $K\in \Cyk$.
\bp\label{prop-cyl-form} For any $K\in \Cyk$, one has
\begin{equation} \label{cha-cyl-1} K=\bigcap_{u\in \sphere} C^{-}(u,
c_K(u)).\end{equation}  Moreover, for any $x\in \partial K$, there
must exist $u_0\in \sphere$ such that
\begin{equation}\label{bounary-sin-1} [x, u_0]=c_K(u_0). \end{equation}  \ep

\begin{proof} It follows from \eqref{cyl-supp-func} that, for a given $u\in \sphere$, $[y, u] \leq c_K(u)$ holds  for any $y\in K$. Thus,
$K\subseteq C^{-}(u, c_K(u))$ and $K\subseteq \bigcap_{u\in \sphere}
C^{-}(u, c_K(u)).$

Assume that $K\subsetneq \bigcap_{u\in \sphere} C^{-}(u, c_K(u)).$
Then there exists $x\in  \bigcap_{u\in \sphere} C^{-}(u, c_K(u))$
but $x\notin K$. Note that $K\in \Cyk$ is a convex body obtained by
the intersection of closed solid cylinders, say $$K=\bigcap_{u\in
\Omega\subseteq \sphere} C^{-}(u, r(u)),$$ where $r(u)$ depends on
$u\in \Omega$. Thus, $x\notin K$ means that there exists $u_0\in
\Omega$ such that $[x, u_0]>r(u_0)$. Note that  $[x, u_0]\leq
c_K(u_0)$ and thus $ r(u_0)<c_K(u_0)$. As $K\subset C^{-}(u_0,
r(u_0))$, then
$$\max_{y\in K}\,[u_0, y]\leq r(u_0)<c_K(u_0),$$ which contradicts
with the maximality of $c_K(u_0)$. Hence \eqref{cha-cyl-1}  follows.

Now let us prove \eqref{bounary-sin-1} for $x\in \partial K$ and
some $u_0\in \sphere$. As $K\in \Cyk$ is an origin-symmetric convex
body, then the convex function $c_K$ restricted on $\sphere$ is
continuous and bounded. Moreover, $x\in \partial K$ implies that
$\lambda x\notin K$ for all $\lambda>1$. Thus, it follows from
\eqref{cha-cyl-1} that there exists $u_{\lambda}\in \sphere$ such
that $ [\lambda x, u_{\lambda}]>c_K(u_{\lambda}).$  As $\sphere$ is
compact and $c_K$ restricted on $\sphere$ is continuous, one can
select a sequence $\lambda_m=1+\frac{1}{m}\rightarrow 1$ whose
corresponding $u_m$ satisfies $u_m\rightarrow u_0$ for some $u_0\in
\sphere$ and $c_K(u_m)\rightarrow c_K(u_0)$ as $m\rightarrow\infty$.
It then follows that
$$[x, u_0] =\lim_{m\rightarrow \infty}[\lambda_mx, u_m]\geq
\lim_{m\rightarrow\infty} c_K(u_m)=c_K(u_0).$$ On the other hand, $[x, u_0]\leq
c_K(u_0)$ as
$x\in K=\bigcap_{u\in \sphere} C^{-}(u, c_K(u)).$ This concludes that if $x\in \partial K$, there must have
$u_0\in \sphere$ such that $[x, u_0]=c_K(u_0)$.  \end{proof}

\bd  Let $K\in \Cyk$  with $n\geq 2$.  For $u\in \sphere$, define the
supporting cylinder of $K$ at direction $u$ by $$  C(u, c_K(u))=\big\{z\in \Rn: [z, u]=c_K(u)\big\}.$$
\ed

Note that  $c_K(u)$ is just the distance from the origin to the
supporting cylinder $C(u, c_K(u))$. Moreover,
 the proof of Proposition \ref{prop-cyl-form} also gives  $c_K(u)=\min\big\{r>0: K\subseteq C^{-}(u, r)\big\}.$
   Hence, $c_K(u)$  is the minimal base radius of the cylinder with axis $\{tu: t\in \R\}$. From \eqref{cha-cyl-1} and its proof,
    if $K\in \Cyk$ and $x\in \partial K$, then there must have (at least) one supporting cylinder, say $C(u, c_K(u))$, containing $x$.
    This observation can be used to define the so-called cylindrical Gauss image of $K$ and its reverse.  More precisely,
   the {\em  cylindrical Gauss image} of $K$,  denoted by $\xi_K: \partial K\rightarrow \sphere$, is defined by, for $x\in \partial K$,   $$\xi_K(x)=\big\{u\in \sphere: u \ \ \mathrm{satisfies} \ \ [x, u] =c_K(u)\big\}.$$ The {\em reverse cylindrical Gauss image} of $K$, denoted by
    $\xi_K^{-1}: \sphere\rightarrow \partial K$, is given by, for $u\in \sphere$,
    $$\xi_K^{-1}(u)=\big\{x\in \partial K:  x \ \ \mathrm{satisfies}\ \  [u, x]=c_K(u)\big\}.$$ Of course, both $\xi_K$ and  $\xi_K^{-1}$ may not be injective. Moreover,  $[u, \xi_K^{-1}(u)]=c_K(u)$ and
    $[x, \xi_K(x)]=c_K(\xi_K(x))$, if $\xi_K(x)$ and $\xi_K^{-1}(u)$ are both singleton sets.

By \eqref{rh} and \eqref{volume-12-5}, the volume of $\polar$ for
$K\in \cKon$ can be calculated by
\begin{equation}\label{povol}V(\polar)=\omega_n\int_{\sphere} \frac{1}{h_K(u)^n}\,du.\end{equation} A similar result also holds for $\ksp$ if $K\in \Cyk$.
\bp\label{prop-1-10-20} Let $K\in \Cyk$  with $n\geq 2$.  Then
$$\rho_{\ksp}(u)\cdot c_K(u)=\rho_{K}(u)\cdot c_{\ksp}(u)=1$$ holds
for all $u\in \sphere$. Moreover,
\begin{equation}\label{volume-sine-polar-12-27} V(\ksp)=\omega_n
\int_{\sphere} \frac{1}{c_K(u)^n}\,du.\end{equation}   \ep

\begin{proof}  By \eqref{cyl-supp-func}, one has, for a given $u\in \sphere$, $[y, u] \leq c_K(u)$  for any $y\in K$.  Equivalently, $$\Big[y, \frac{u}{c_K(u)} \Big]\leq 1\ \ \mathrm{for\  all }\ y\in K.$$ Hence $ \frac{u}{c_K(u)} \in \ksp$ by Definition \ref{s-p-b}. It follows from \eqref{radial} that $\rho_{\ksp}(u)\geq \frac{1}{c_K(u)}$. On the other hand, one can take $y_0\in \xi_K^{-1}(u)$. Therefore, $[y_0,  \frac{u}{c_K(u)} ]=1$ and $[y_0,  tu ]>1$ if $t>\frac{1}{c_K(u)}$. As $y_0\in \partial K$, it follows from Definition \ref{s-p-b} that $tu\notin \ksp$ for any $t>\frac{1}{c_K(u)}$. The fact that $\rho_{\ksp}(u)u\in \partial \ksp$ then  yields  $\rho_{\ksp}(u)\leq \frac{1}{c_K(u)}$  and thus $\rho_{\ksp}(u)= \frac{1}{c_K(u)}$. By Proposition \ref{prop-p-s} (v), one gets that, for $u\in \sphere$, $$\rho_K(u)=\rho_{K^{\diamond\diamond}}(u)= \frac{1}{c_{\ksp}(u)}.$$ Finally, the volume of $\ksp$ can be calculated by $$V(\ksp)=\omega_n\int_{\sphere} \rho_{\ksp}(u)^n\,du=\omega_n\int_{\sphere}   \frac{1}{c_K(u)^n} \,du$$
 due to \eqref{volume-12-5}.\end{proof}

After the above preparations, we now prove Proposition
\ref{prop-p-s} (vi).
\begin{proof}[Proof of  Proposition \ref{prop-p-s} (vi)]
        Recall that $K=\bigcap_{u\in \sphere} C^{-}(u, c_K(u))$ proved in
    \eqref{cha-cyl-1}.   Then $K\subset C^{-}(u, c_K(u))$ for any fixed
    $u\in \sphere$ and  $[x, u]\leq c_K(u)$ for any $x\in K$. By
    \eqref{cal-sin}, one must have
    \begin{equation*}[x, u]= |\textrm{P}_{u^{\perp}}x|\geq |\textrm{P}_v x|
        =|x\cdot v|\end{equation*} for any $v\in \sphere$ such that $u\perp
    v$. Taking the maximum over $x\in K$, one has
    \begin{equation*}
        c_K(u)\geq h_K(v),\end{equation*}
    and thus,  for any $v\in S^{n-1}$,
    \begin{equation}\label{compare-two-support}
    \int_{S^{n-1}\cap v^\perp}\frac{1}{c_K(\zeta)^n}\,d\zeta\leq\frac{1}{h_K(v)^n},
    \end{equation}
    where $d\zeta$ is the rotation invariant probability measure on $S^{n-1}\cap v^\perp$. For any nonnegative continuous function  $f:
    \sphere\rightarrow \R$, \cite[(2.22)]{Koldobsky-FA0} yields that
    \begin{equation}\label{int-grass-1} \int_{\sphere}
        f(u)\,du=\int_{S^{n-1}}\bigg(\int_{S^{n-1}\cap v^\perp} f(\zeta)\, d\zeta\bigg)\, dv.\end{equation} By
    \eqref{volume-sine-polar-12-27}, \eqref{int-grass-1}, \eqref{compare-two-support}, and \eqref{povol}, one
    gets
    \begin{eqnarray*} V(K^{\diamond}) \!\!&=&\!\!\omega_n \int_{\sphere}
        \frac{1}{c_K(u)^n}\,du \nonumber
        \\\!\!&=&\!\!\omega_n \int_{S^{n-1}}
        \bigg(\int_{\sphere\cap v^{\perp}}
        \frac{1}{c_K(\zeta)^n} \,d\zeta\bigg)\,dv \nonumber \\
        \!\!&\leq &\!\!\omega_n \int_{S^{n-1}}
        \frac{1}{h_K(v)^n}\,dv,  \nonumber \\
    \!\!&= &\!\! V(K^{\circ}).\end{eqnarray*}
     This concludes the proof of
\eqref{compare-polar-diamond}.
    \end{proof}

  \section{The $L_p$-sine Blaschke-Santal\'{o} inequality}\label{sec-3-2}

In this section, we will prove the $L_p$-sine
Blaschke-Santal\'{o} inequality \eqref{main2}. The following result is needed. 

\bp\label{calculate-the-ball}  Let $p\ge 1$ and $n\ge 2$. Then
\begin{equation}\label{B}
    \Lambda_p \ball=\ball.
\end{equation}
\ep

\begin{proof}  Lemma 7.5 in \cite{LXZ} states that for an $m$-dimensional subspace  $V$  of $\mathbb R^n$, one has
\begin{eqnarray}\label{sphere}
\int_{\sphere}\big|\textrm{P}_Vu\big|^p\,du=\frac{m\omega_m\omega_{n+p-2}}{n\omega_n\omega_{m+p-2}}
\ \ \mathrm{for} \ \ p\geq 1,
\end{eqnarray}
where $\textrm{P}_{V}u$ is the orthogonal projection of $u$ onto the
subspace $V$. Recall that $[x, u]=|x|\cdot
|\textrm{P}_{x^{\perp}}u|$ for all $x\in \Rn\setminus\{o\}$ (see
\eqref{cal-sin}). Hence, if $x\in \Rn\setminus\{o\}$, one can let
$V=x^{\perp}$ whose dimension is $m=n-1$ and then \eqref{sphere}
yields  that for all $p\geq 1$,
\begin{equation*}
\int_{\sphere}[x, u]^p\,du=|x|^p\cdot
\int_{\sphere}|\textrm{P}_{x^{\perp}}u|^p\,du=\bigg(\frac{(n-1)\omega_{n-1}\omega_{n+p-2}}{n\omega_n\omega_{n+p-3}}\bigg)
|x|^p.
\end{equation*}
Together with \eqref{def-p-sine-1} and \eqref{tcnp},  one gets
\begin{eqnarray*}
    h_{\Lambda_p\ball}(x)\!\!&=&\!\!\Big(\frac{n\omega_n}{(n+p)\widetilde{c}_{n,p}V(B^n)}\int_{\sphere}[x,
    u]^p\rho_{\ball}(u)^{n+p}\,du\Big)^{1/p}
    \\\!\!&=&\!\!\Big(\frac{n}{(n+p)\widetilde{c}_{n,p}}\int_{\sphere}[x, u]^p\,du\Big)^{1/p}\\\!\!&=&\!\!|x|=h_{\ball}(x).  \end{eqnarray*} This concludes \eqref{B}.
     \end{proof}
Moreover, it follows from \eqref{c} and \eqref{B} that
\begin{equation}\label{cB}
    \Lambda_p^{\circ} (c\ball)=c^{-1}\Lambda_p^{\circ} (\ball)=c^{-1}(\ball)^{\circ}=c^{-1}\ball,\ \ \ \mathrm{for \ any}\ c>0.
\end{equation}

\bl\label{LL} Let $p\ge 1$ and $n\ge 2$. If $K\subset\Rn$ is an $L_{n+p}$-star  with $V(K)>0$, then
\begin{equation}\label{L2-12-12}
    V(\Lambda_p^{\circ}\plp)\geq V(K),
\end{equation}
with equality if and only if
$\rho_K(u)= \rho_{\Lambda_p^{\circ}\plp}(u)$ for almost all $u\in \sphere$ with respect to the spherical measure on $\sphere$.

\el

 \begin{proof} Let $K, L\subset\Rn$ be  two $L_{n+p}$-stars  with positive volumes. As discussed  in Section \ref{sec-3-1}, both $\plp\in \cKe$ and $\Lambda_p^{\circ}L\in \cKe$ are origin-symmetric convex bodies. It follows from \eqref{dual}, \eqref{Po}, and Fubini's theorem that
\begin{eqnarray}
\frac{\widetilde{V}_{-p}(K,\Lambda_p^{\circ}L)}{V(K)}\!\!&=&\!\!\frac{\omega_n}{V(K)}\int_{\sphere}\rho_K(u)^{n+p}\rho_{\Lambda_p^{\circ}L}(u)^{-p}\,du
\nonumber
\\\!\!&=&\!\!\frac{n\omega_n^2}{(n+p)\widetilde{c}_{n,p}V(K)V(L)}\int_{\sphere}\int_{\sphere}[u,v]^p\rho_K(u)^{n+p}\rho_{L}(v)^{n+p}\,dv\,du
\nonumber \\\!\!&=&\!\!
\frac{\omega_n}{V(L)}\int_{\sphere}\rho_{L}(v)^{n+p}\Big(\frac{n\omega_n}{(n+p)\widetilde{c}_{n,p}V(K)}\int_{\sphere}[u,v]^p
\rho_K(u)^{n+p}\,du\Big)\,dv\nonumber
 \\\!\!&=&\!\!\frac{\omega_n}{V(L)}\int_{S^{n-1}}\rho_L(v)^{n+p}\rho_{\plp}(v)^{-p}\,dv \nonumber
   \\\!\!&=&\!\!\frac{\widetilde{V}_{-p}(L,\plp)}{V(L)}.\label{rel}
\end{eqnarray} By letting $L=\plp$,   \eqref{e-21} and \eqref{rel} imply that
\begin{eqnarray}
\frac{\widetilde{V}_{-p}(K,\Lambda_p^{\circ}\plp)}{V(K)}&=&\frac{\widetilde{V}_{-p}(\plp,\plp)}{V(\plp)}=1.\label{rel21}
\end{eqnarray} That is, $\widetilde{V}_{-p}(K,\Lambda_p^{\circ}\plp)=V(K).$ Together with  \eqref{dm}, one has, for all $p\ge 1$,
\begin{equation*} V(K)^n
= \widetilde{V}_{-p}(K, \Lambda_p^{\circ}\plp)^n \geq   V(K)^{n+p}V(\Lambda_p^{\circ}\plp)^{-p},
\end{equation*} which is exactly \eqref{L2-12-12} after rearrangement. Clearly if $\rho_K(u)= \rho_{\Lambda_p^{\circ}\plp}(u)$ for almost all $u\in\sphere$ with respect to the spherical measure of $\sphere$, then  $V(\Lambda_p^{\circ}\plp)= V(K)$. On the other hand, if
  $V(\Lambda_p^{\circ}\plp)= V(K)$ is assumed, then the equality conditions of \eqref{dm} yield that there exists a constant $c>0$ such that  $\rho_K(u)=c\cdot \rho_{\Lambda_p^{\circ}\plp}(u)$ for almost all $u\in \sphere$ with respect to the spherical measure of $\sphere$. This further yields that $V(\Lambda_p^{\circ}\plp)=V(K)=c^nV(\Lambda_p^{\circ}\plp)$ and hence $c=1$, as desired.
\end{proof}

Denote by $[x_1,\cdots, x_n]$ the $n$-dimensional volume of the
parallelotope spanned by vectors $x_1, \cdots, x_n\in \Rn$. Let
$p\ge 1$ and $n\ge 2$.  We shall need $\mathrm{T}_p(K_2,\cdots,
K_n)\in \cKe$ for $K_2,\cdots, K_n\in \cKon$, whose support function
at $x\in \Rn$ was given in \cite{LHX} (up to a factor) by
\begin{equation*}
    h_{\mathrm{T}_p(K_2,\cdots, K_n)}(x)^p=\frac{1}{V(K_2)\cdots V(K_n)}\int_{K_2}\cdots\int_{K_n}[x,x_2,\cdots,x_n]^p\,dx_2\cdots
    \,dx_n.
\end{equation*}
In particular, if we let $K_2=K\in \cKon$, $K_3=\cdots =K_n=\ball$,
then \cite[Theorem 4.3]{LHX} asserts that
\begin{equation*}\mathrm{T}_p(K, \underbrace{\ball, \cdots,
\ball}_{n-2})= d_{n, p}\cdot\Lambda_{p}K
\end{equation*} for all $p\geq 1$ with $d_{n, p}^{p}= \widetilde{c}_{n,p}$ for $n=2$ and
 \begin{eqnarray}  \frac{(n+p)^{n-2}d_{n, p}^p}{n^{n-2} \widetilde{c}_{n,p}}\!\!&=&\!\!\int_{\sphere} \cdots \int_{\sphere} \bigg(\frac{[x, u_2, \cdots, u_n]}{[x, u_2]}\bigg)^p\,du_2\cdots \,du_n\nonumber  \\\!\!&=&\!\!
   \prod_{i=2}^{n-1}  \int_{\sphere}\big|\textrm{P}_{V_i^{\perp}}u_{i+1}\big|^p\,du_{i+1}= \prod _{i=2}^{n-1} \frac{(n-i)\omega_{n-i}\omega_{n+p-2}}{n \omega_n\omega_{n-i+p-2}}\label{formular419} \end{eqnarray}
for $n\geq 3$. The second equality in \eqref{formular419} follows
from \cite[(3.7)]{LHX}, where the dimensions of subspaces
$V_i^{\perp}$ are $n-i$ for $i=2, \cdots, n$. (We would like to
point out that $\,dv_i$ in
   \cite[Lemma 3.2]{LHX}  are equal to  $du_i$ in the present paper multiplying  $n\omega_n$.) The last equality in \eqref{formular419} follows
from formula \eqref{sphere}. Therefore, by \eqref{tcnp}, one  has
\begin{equation*} d_{n, p}  =  \bigg(\prod _{i=1}^{n-1}\frac{(n-i)\omega_{n-i}\omega_{n+p-2}}{(n+p)\omega_n\omega_{n-i+p-2}}\bigg)^{1/p}.\end{equation*}
By \eqref{cm}, we
further have, for $p\geq 1$ and $n\geq 2$,
\begin{eqnarray}
    \mathrm{T}_{p}^{\circ}(K, \underbrace{\ball,\cdots,\ball}_{n-2})=d_{n, p}^{-1}\cdot \Lambda_{p}^{\circ}K. \label{relation-2}
   \end{eqnarray}

With a different normalization, the following inequality was
established by Haddad in \cite[Theorem 1.3]{H}, which is equivalent
to the $L_p$-Busemann random simplex inequality \cite[Theorem
1.3]{BMMP}. By a different approach, inequality \eqref{cen1} for
star bodies and for $p=2$ was also established in \cite[Corollary
4.5]{LHX}.

\bl \label{t-1} Let $p\geq 1$ and $n\geq 2$. If $K_2,\cdots,K_n\in
\cKon$, then
\begin{equation}\label{cen1}
V(K_2)\cdots V(K_n)\cdot V(\mathrm{T}_p^{\circ}(K_2,\cdots,
K_n))\leq
\bigg(\frac{n\omega_n^{n+p}}{(n+p)\omega_{n+p}^n}\Big(\prod_{j=1}^{n}\frac{(p+j)\omega_{p+j}}{j\omega_j}\Big)\bigg)^{\frac{n}{p}},
\end{equation} with equality if and only if $K_2,\cdots,K_n$ are
origin-symmetric ellipsoids that are dilates. \el

Before the proof of Theorem \ref{lpm}, we shall mention the  $L_p$-sine centroid body and the sine polar body for the case $n=2$. Notice that  for $x,y\in\mathbb{R}^2$,
\begin{equation}\label{translation-90} [x, y]=|x|\cdot |\textrm{P}_{x^{\perp}}y|=|\psi_{\pi/2}x|\cdot |\textrm{P}_{\psi_{\pi/2}x}y|=|\psi_{\pi/2}x\cdot y|,\end{equation}
where 
\begin{equation*}\psi_{\pi/2}=\left(
	\begin{array}{ccc}
		0 & -1\\
		1& 0 \\
	\end{array}
	\right)\end{equation*}  is a rotation by the angle of~$\pi/2$  in
$\mathbb{R}^2$. 
For $K\in \mathscr{K}_{o}^2$, it follows from \eqref{cnp} and
\eqref{tcnp} that $c_{2,p}=\widetilde{c}_{2,p}$. By
\eqref{def-p-sine-1}, \eqref{translation-90}, \eqref{centroid}, and
\eqref{affh}, we have, for $p\geq 1$  and for $x\in \mathbb R^2$,
\begin{equation*}  h_{\Lambda_pK}(x)^p=\frac{1}{\widetilde{c}_{2,p}V(K)}\int_K[x,y]^p\,dy  =\frac{1}{c_{2,p}V(K)}\int_K |\psi_{\pi/2}x\cdot
	y|^p\,dy=h_{\psi_{\pi/2}^{t}\Gamma_pK}(x)^p. \end{equation*}
Together with \eqref{19} and the the facts that
$\psi_{\pi/2}^{-1}=-\psi_{\pi/2}$ and $\Gamma_p^{\circ}K$ is
origin-symmetric, one has
\begin{equation}\label{e-2}
	\Lambda_{p}^{\circ}K=\psi_{\pi/2}\Gamma_p^{\circ}K.
\end{equation}
Similar, for $K\in \mathscr{K}_{e}^2$, it follows from \eqref{sc},
\eqref{polar}, and \eqref{translation-90} that
\begin{equation}\label{rot}
	K^{\diamond}=\psi_{\pi/2} \polar.
\end{equation}
Due to the affine natures of $\Gamma_p^{\circ}K$ and $\polar$, it
is  easy to verify that in $\mathbb{R}^2$
\begin{equation}\label{lambaff}
	\Lambda_{p}^{\circ}(\phi K)=|\det\phi|^{-1}\phi \Lambda_{p}^{\circ}K,
\end{equation}
and
\begin{equation}\label{poaff}
	(\phi K)^{\diamond}=|\det\phi|^{-1}\phi K^{\diamond}.
\end{equation}
Therefore, the volume product~$V(K)V(\Lambda_{p}^{\circ}K)$
is~$\textrm{GL}(2)$-invariant for $K\in \mathscr{K}_{o}^2$ and the
volume product~$V(K)V(\ksp)$ is~$\textrm{GL}(2)$-invariant  for
$K\in \mathscr{K}_{e}^2$. The volume products
$V(K)V(\Lambda_{p}^{\circ}K)$ and~$V(K)V(\ksp)$, however, is no
longer~$\textrm{GL}(n)$-invariant for $n\geq 3$. This is a major
difference between equality conditions in inequality \eqref{pBS},
\eqref{BS-polar-1} and its sine counterpart \eqref{main2},
\eqref{SB}.

We are now in a position to prove the $L_p$-sine
Blaschke-Santal\'{o} inequality \eqref{main2} in the following
theorem. For $p=2$ and $K$ being a star body, one recovers
\cite[Theorem 4.1]{LHX}.

\bt\label{lpm} Let $p\ge 1$ and $n\ge 2$. If $K\subset\Rn$ is an
$L_{n+p}$-star  with $V(K)>0$, then
\begin{equation}\label{sine}
V(K)V(\plp)\leq \omega_n^{2},
\end{equation} with equality if and only if $K$, up to sets of measure $0$, is an
origin-symmetric ellipsoid when $n=2$, and is an origin-symmetric
ball when $n\geq3$.
 \et

\begin{proof}
Note that when $n=2$ and $K\in \mathscr{K}_{o}^2$, it follows from 
\eqref{e-2} that inequality \eqref{sine} is exactly the $L_p$
Blaschke-Santal\'{o} inequality \eqref{pBS}.

Now assume that $n\geq 3$ and $K\in \cKon$. Taking $K_2=K$ and
$K_3=\cdots=K_n=\ball$ in \eqref{cen1}, one gets
\begin{equation*} V(K)\cdot \underbrace{V(\ball)\cdots
V(\ball)}_{n-2} \cdot V(\mathrm{T}_p^{\circ}(K, \underbrace{\ball
\cdots, \ball}_{n-2}))\leq
\bigg(\frac{n\omega_n^{n+p}}{(n+p)\omega_{n+p}^n}\Big(\prod_{j=1}^{n}\frac{(p+j)\omega_{p+j}}{j\omega_j}\Big)\bigg)^{\frac{n}{p}},
\label{sine-12-27-0}
\end{equation*} with equality if and only if $K_2=K$ and $K_3=\cdots=K_n=\ball$ are
origin-symmetric ellipsoids that are dilates (of course, this
implies that $K$ is an origin-symmetric ball).

 Together with \eqref{relation-2}, one has, for $n\geq 3$,
\begin{eqnarray*}
V(K)V(\plp) \!\!&\leq&\!\! \frac{d_{n,
p}^{n}}{\omega_n^{n-2}}\bigg(\frac{n\omega_n^{n+p}}{(n+p)\omega_{n+p}^n}\Big(\prod_{j=1}^{n}\frac{(p+j)\omega_{p+j}}{j\omega_j}\Big)\bigg)^{\frac{n}{p}}
\nonumber \\\!\!&=&\!\!\omega_n^2 \bigg(\frac{(p+1)\cdots (p+n)
\omega_{n+p-1}\omega_{n+p-2}^{n}}{
(n+p)^{n}\omega_p\omega_{p-1}\omega_{n+p}^{n-1} }
\bigg)^{\frac{n}{p}} \nonumber \\\!\!&=&\!\!\omega_n^{2},
\end{eqnarray*} with equality if and only if $K$ is an
origin-symmetric ball.

 In the above, we have established inequality \eqref{sine} for
all $n\geq 2$ and $K\in \cKon$. Now we consider $K$ being an
$L_{n+p}$-star  with $V(K)>0$. Together with \eqref{L2-12-12} and
the fact that $\plp\in \cKon$, we have \begin{equation}\label{gee}
  V(K)V(\plp)\le V(\Lambda_{p}^{\circ}\Lambda_{p}^{\circ}K)V(\plp)\leq  \omega_n^{2}.
\end{equation} Hence, the desired inequality \eqref{sine}  holds for $n\geq 2$ and for any $L_{n+p}$-star  $K$ with $V(K)>0$.

 Now let us characterize the equality of \eqref{sine}.
Equality holding in the second inequality of \eqref{gee} yields that
$\plp$ is an origin-symmetric ellipsoid when $n=2$ and is an
origin-symmetric ball when $n\geq3$. For $n=2$, we may assume
$\plp=\phi B^2$  for some $\phi\in \mathrm{GL}(2)$. By
\eqref{lambaff} and \eqref{cB}, one has
\begin{equation}\label{lam2}
	\Lambda_{p}^{\circ}\Lambda_{p}^{\circ}K =\Lambda_{p}^{\circ}(\phi
	B^2)=|\det \phi|^{-1}\phi\Lambda_p^{\circ}B^2=|\det \phi|^{-1} \phi
	B^2.\end{equation} For $n\ge 3$, we may assume $\plp=c\ball$ for some
constant $c>0$.   By \eqref{cB}, one has
\begin{equation*}
	\Lambda_{p}^{\circ}\Lambda_{p}^{\circ}K=\Lambda_{p}^{\circ}(c\ball)=c^{-1}
	\ball. \label{bal-12-27} \end{equation*}  From Lemma \ref{LL},
equality holding in the first inequality of \eqref{gee} implies
$\rho_K(u)= \rho_{\Lambda_p^{\circ}\plp}(u)$ for almost all $u\in
\sphere$ with respect to the spherical measure on $\sphere$.
That is, to have equality  in \eqref{gee},  $K$, up to sets of measure $0$, is an
origin-symmetric ellipsoid when $n=2$, and is an origin-symmetric
ball when $n\geq3$. Conversely, due to \eqref{e-2} and \eqref{cB},
equality holds in \eqref{sine} if  $K$, up to sets of measure $0$,
is an origin-symmetric ellipsoid when $n=2$, and is an
origin-symmetric ball when $n\geq3$.
\end{proof}

\section{The sine Blaschke-Santal\'{o} inequality} \label{BS-Sine-2-2-1}

If $K$ is a star body in~$\Rn$, then $K$ is an $L_{n+p}$-star for
all $p\geq 1$ with $V(K)>0$. By taking $p\rightarrow \infty$ on
\eqref{sine} and \eqref{dinf}, one gets the sine Blaschke-Santal\'{o}
inequality \eqref{SB}:
\begin{equation}\label{Lambda-Infty-1} V(K)V(\ksp)=V(K)V(\pli) =  \lim_{p\rightarrow \infty} V(K)V(\plp)  \leq  \omega_n^{2}.
\end{equation} However, this approximation
argument does not yield equality conditions of
\eqref{Lambda-Infty-1}. The full characterization of equality
conditions of \eqref{Lambda-Infty-1} will be presented in Theorem
\ref{santalo-sine-12-27}.

The following lemma will be used to extend the sine
Blaschke-Santal\'{o} inequality \eqref{Lambda-Infty-1}  from star
bodies  to general bounded and measurable sets in $\Rn$.

\bl\label{dd} If $K$ is a bounded and measurable set in $\Rn$ with
$n\geq 2$, then
\begin{equation}\label{L2-infty}
   V(K^{\diamond\diamond})\geq V(K).
\end{equation}
If in addition $K$ is a star body, equality holds in
\eqref{L2-infty} if and only if $K=K^{\diamond\diamond}$. \el

 \begin{proof} As $K$ is bounded, $K^{\diamond\diamond}$ is also bounded by Proposition \ref {prop-p-s}  (i)-(iii).  Proposition \ref {prop-p-s} (v) implies $K\subseteq
K^{\diamond\diamond}$ and hence $V(K^{\diamond\diamond})\geq V(K)$.
Now suppose  $V(K^{\diamond\diamond})=V(K)$ and  $K$ is a star body
in $\Rn$. It follows from \eqref{volume-12-5} that
$$0=V(K^{\diamond\diamond})-V(K)=\omega_n \int _{\sphere}
\Big(\rho_{K^{\diamond\diamond}}(u)^n-\rho_K(u)^n\Big)\,du.$$ Since
$\rho_{K^{\diamond\diamond}}\geq \rho_K$ and $\rho_K$ is continuous
on $\sphere$, this further implies that
$\rho_{K^{\diamond\diamond}}(u)= \rho_K(u)$ for all $u\in \sphere$,
and thus $K^{\diamond\diamond}=K$. 
\end{proof}

We are now in a position to prove the sine Blaschke-Santal\'{o}
inequality \eqref{SB}.

 \bt \label{santalo-sine-12-27} Let $K$ be a bounded and measurable set in $\Rn$ with $n\geq 2$ such that $K$ is not concentrated on any $1$-dimensional subspace. Then
 \begin{equation}\label{sine-1}
V(K)V(\ksp) \leq   \omega_n^{2}.
\end{equation} If in addition $K$ is a star body, equality holds in \eqref{sine-1} if and only if $K$  is an
origin-symmetric ellipsoid when $n=2$ and is an origin-symmetric ball when $n\geq3$. \et

 \begin{proof} Since $K\subset \Rn$ is bounded and is not concentrated on any $1$-dimensional subspace, it follows from Proposition \ref{prop-p-s} (iv) that $K^{\diamond}\in \Cyk$ and $V(K^{\diamond})<\infty$. Clearly, inequality \eqref{sine-1} holds trivially
 if~$V(K)=0$.

Assume that $V(K)>0$. Lemma \ref{dd} yields
\begin{equation}\label{V2}V(K) V(\ksp) \leq
	V(\ksp)V(K^{\diamond\diamond}).\end{equation} If
$K\in\mathcal{S}^n$, then equality in \eqref{V2} holds if and only
if $K^{\diamond\diamond}=K$. Applying \eqref{compare-polar-diamond}
for~$L=K^{\diamond}\in \Cyk$~and the Blaschke-Santal\'{o} inequality
\eqref{BS-polar-1} for~$\ksp$, one gets
\begin{equation}\label{V31}V(\ksp)V(K^{\diamond\diamond})\leq V(\ksp)V((K^{\diamond})^{\circ})\leq \omega_n^2.\end{equation}
Thus, the desired inequality \eqref{sine-1} immediately follows from
\eqref{V2} and \eqref{V31}.

Now let us characterize the equality of \eqref{sine-1} under the
assumption that $K$ is a star body.  Clearly, from \eqref{rot},
Proposition \ref {prop-p-s} (ii) and (iii), equality holds in
\eqref{sine-1} if  $K$ is an origin-symmetric ellipsoid when $n=2$
and is an origin-symmetric ball when $n\geq3$. Conversely, the
equality conditions of the Blaschke-Santal\'{o} inequality
\eqref{BS-polar-1} imply $\ksp$ is an origin-symmetric ellipsoid.  
When $n=2$, we may assume $\ksp=\phi B^2$ for some $\phi\in
\mathrm{GL}(2)$. It follows from the equality conditions of
\eqref{V2}, \eqref{poaff} and Proposition \ref{prop-p-s} (ii) that
\begin{equation*}K=K^{\diamond\diamond}=(\phi B^2)^{\diamond}=(\det\phi)^{-1}\phi (B^2)^{\diamond}=(\det\phi)^{-1}\phi B^2,\end{equation*}
which is an origin-symmetric ellipsoid. To get the desired equality conditions for $n\geq 3$, we
only need to show that equality holding in \eqref{sine-1} for~$n\ge
3$ yields that $\ksp$ is an origin-symmetric ball. Suppose
that~$K^{\diamond}$~is  an
origin-symmetric ellipsoid but is not an origin-symmetric ball. Without loss of
generality, one can assume that
\begin{equation*}\label{ell}\ksp=\mathscr{E}=\Big\{x=(x_1, \cdots, x_n)\in
\Rn: \ \
\frac{x_1^2}{a_1^2}+\frac{x_2^2}{a_2^2}+\cdots+\frac{x_n^2}{a_n^2}\le
1\Big\},
\end{equation*} where $a_1\le a_2\le \cdots\le a_n$ such that  $a_1<a_n$.
We separate the proof into two cases.

\vskip 2mm \noindent {\it Case 1:}  there exists $1< k<n$ such that
$a_1< a_k\leq a_{n}$.  In this case, we shall prove that
$\mathscr{E}\notin \Cyk$, which is a contradiction with the fact
$\mathscr{E}=\ksp\in \Cyk$.

We argue by contradiction by assuming that $\mathscr{E}\in \Cyk$.
Note that $a_1e_1\in \partial \mathscr{E}$. It follows from
Proposition \ref{prop-cyl-form}  that  there is some $u_0\in\sphere$
such that
\begin{equation}\label{ellip-u-0}  c_{\mathscr{E}}(u_0)=[a_1e_1,
u_0]=a_1[e_1, u_0]\leq a_1.\end{equation} Let $z=(z_1, \cdots,
z_n)\in u_0^{\perp}\cap
\partial \mathscr{E}$ such that $z\notin \textrm{span}\{e_1,\ldots,e_{k-1}\}$, where $\{e_1, \cdots, e_n\}$ is the canonical basis of $\Rn$. Since $z\notin  \textrm{span}\{e_1,\ldots,e_{k-1}\}$, then at least one of $z_{k}, \cdots, z_n$ is not zero. Together with $z\in \partial \mathscr{E}$ and $a_1< a_k\leq a_{n}$, one can check that $$ 1=\frac{z_1^2}{a_1^2}+\cdots+\frac{z_k^2}{a_k^2}+\cdots+\frac{z_n^2}{a_n^2}<\frac{|z|^2}{a_1^2},
$$  which implies $|z|>a_1$. On the other hand,  due to $z\in u_0^{\perp}$,  one gets $[z, u_0]=|z|>a_1.$ It follows from \eqref{cyl-supp-func}  that
\begin{equation*} c_{\mathscr{E}}(u_0)=\max_{y\in \mathscr{E}} [y,  u_0]\geq [z, u_0]>a_1.\end{equation*} This contradicts to \eqref{ellip-u-0} and hence $\mathscr{E}\notin\Cyk$.

 \vskip 2mm \noindent {\it Case 2:}  let $a:=a_1=\cdots=a_{n-1}$ and $b:=a_n$ with $a<b$. Then \begin{equation} \mathscr{E}=\Big\{x=(x_1, \cdots, x_n)\in \Rn: \ \ \frac{x_1^2}{a^2}+\cdots+\frac{x_{n-1}^2}{a^2}+\frac{x_n^2}{b^2}\le
1\Big\}. \label{spe-ellp-1}\end{equation} It is well-known that
\begin{equation}\label{Ep}\mathscr{E}^{\circ}=\Big\{x=(x_1, \cdots, x_n)\in \Rn: \ \ a^2x_1^2+\cdots+a^2x_{n-1}^2+b^2x_n^2\le
1\Big\}. \end{equation} We now claim that \begin{equation}\label{Ed}\mathscr{E}^{\diamond}\subseteq \Big\{x=(x_1, \cdots, x_n)\in \Rn: \ \ b^2x_1^2+\cdots+b^2x_{n-1}^2+a^2x_n^2\le
1\Big\}:=\mathscr{E}_0.
\end{equation}

Let $U_0= \textrm{span}\{e_{n-1}, e_n\}$. Then it follows from
Definition \ref{s-p-b} that
\begin{eqnarray}
   \mathscr{E}^{\diamond}\cap U_0\!\!&=&\!\!\big\{x\in\R^n\cap U_0: \max _{y\in \mathscr{E}}[x, y]\leq 1\big\} \nonumber
   \\\!\!&\subseteq&\!\!\big\{x\in\R^n\cap U_0: \max _{y\in \mathscr{E}\cap U_0}[x, y]\leq 1\big\}. \label{get-ellip-1}
\end{eqnarray} In the $2$-dimensional subspace $U_0$, $\mathscr{E}^{\circ}\cap U_0$ is an ellipse of the following form:
$$\mathscr{E}^{\circ}\cap U_0=\Big\{z=(z_1, \cdots, z_n)\in \Rn: \ \ a^2 z_{n-1}^2+b^2  z_n^2 \le
1, \ \ z_1=\cdots =z_{n-2}=0\Big\}.$$ As shown in \eqref{rot},
restricted on the $2$-dimensional subspace  $U_0$,  the set
$$\mathcal{E}=\big\{x\in\R^n\cap U_0: \max\limits_{y\in
\mathscr{E}\cap U_0}[x, y]\leq 1\big\}$$  is just a rotation of
$\mathscr{E}^{\circ}\cap U_0$ by an angle of $\pi/2$, namely
$$\mathcal{E}=\Big\{z=(z_1, \cdots, z_n)\in \Rn: \ \ b^2
z_{n-1}^2+a^2  z_n^2 \le 1, \ \ z_1=\cdots
=z_{n-2}=0\Big\}=\mathscr{E}_0\cap U_0.$$ Together with
\eqref{get-ellip-1}, one concludes
\begin{equation}\label{compare-spe-ellp-1}
\mathscr{E}^{\diamond}\cap U_0\subseteq
\mathcal{E}=\mathscr{E}_0\cap U_0.\end{equation}  Let $O_n$ be a
rotation of the form
\begin{equation}\label{rotation-n-n-1} O_n=\left(
               \begin{array}{ccc}
                 O_{n-1} & 0\\
                 0& 1 \\
               \end{array}
             \right)\end{equation} where $O_{n-1}$ is an arbitrary rotation on the $(n-1)$-dimensional subspace $\textrm{span}\{e_1, \cdots, e_{n-1}\}.$
             It can be checked by Proposition \ref{prop-p-s} (iii) that $$O_n(\mathscr{E}^{\diamond}\cap U_0) =O_n(\mathscr{E}^{\diamond})\cap O_nU_0=
             (O_n\mathscr{E})^{\diamond}\cap O_nU_0=  \mathscr{E}^{\diamond}\cap O_nU_0,$$ where the last equality follows from the definition of $\mathscr{E}$ defined by \eqref{spe-ellp-1}. Similarly,
 $$O_n(\mathscr{E}_0\cap U_0)=O_n\mathscr{E}_0\cap O_nU_0=\mathscr{E}_0\cap O_nU_0,$$  where the the last equality follows from the definition of $\mathscr{E}_0$ defined by \eqref{Ed}.  Applying $O_n$ on both sides of  \eqref{compare-spe-ellp-1}, one gets
$$\mathscr{E}^{\diamond}\cap O_nU_0=O_n(\mathscr{E}^{\diamond}\cap U_0) \subseteq  O_n(\mathscr{E}_0\cap U_0)=\mathscr{E}_0\cap O_nU_0.$$
Taking the union over all $O_n$ of the form \eqref{rotation-n-n-1}, one gets
$$\mathscr{E}^{\diamond}=\mathscr{E}^{\diamond}\cap\Rn=
\mathscr{E}^{\diamond} \bigcap \Big(\bigcup_{O_n} \big(O_nU_0\big)
\Big)=\bigcup_{O_n} \Big(\mathscr{E}^{\diamond}\cap
O_nU_0\Big)\subseteq \bigcup_{O_n} \Big(\mathscr{E}_0 \cap
O_nU_0\Big)=\mathscr{E}_0,$$ which concludes the proof of
\eqref{Ed}.   Thus, it follows from \eqref{Ep} and \eqref{Ed} that
\begin{equation*}V(\mathscr{E}^{\diamond})\le V(\mathscr{E}_0)=\frac{1}{b^{n-1}a}\omega_n\quad \textrm{and} \quad V(\mathscr{E}^{\circ})=\frac{1}{a^{n-1}b}\omega_n.
\end{equation*}
Since $n\geq 3$ and $b>a$,  one clearly has
$V(\mathscr{E}^{\diamond})<V(\mathscr{E}^{\circ})$. It
means~$V(K^{\diamond\diamond})<
V((K^{\diamond})^{\circ})$~since~$\ksp=\mathscr{E}$. Thus, from
\eqref{V2} and \eqref{V31}, equality cannot hold in inequality
\eqref{sine-1}.
\end{proof}

 \section{Related functional inequalities} \label{section-5}
As applications of  the $L_p$-sine
Blaschke-Santal\'{o} inequality \eqref{sine}, we will provide some  functional inequalities related to it.  These functional inequalities can be proved along the routine approach in the literature, so we only give a detailed proof for Theorem \ref{m1} and omit the proofs for Theorems \ref{m12} and \ref{m2}. 

The following inequality is equivalent to the $L_p$-sine
Blaschke-Santal\'{o} inequality \eqref{sine}.

\bt \label{m1} Let $p\ge 1$ and $n\ge 2$. If  $K$ and $L$ are
$L_{n+p}$-stars  in
$\Rn$, then
\begin{equation}\label{gm}
    \int_K\int_L[x,y]^p\,dx\,dy\geq
    \frac{ n(n-1) \omega_{n-1}\omega_{n+p-2}}{ (n+p)^2 \omega_{n+p-3}\omega_n^{1+2p/n}} \big[V(K)V(L)\big]^{\frac{n+p}{n}}.
\end{equation} If $V(K)V(L)>0$,  equality  holds in \eqref{gm} if and only if $K$ and $L$, up to sets of
measure $0$, are dilates of an origin-symmetric ellipsoid when $n=2$
and are origin-symmetric balls when $n\geq3$. \et

\begin{proof}
    \vskip 2mm Without loss of generality, we assume that $V(K)V(L)>0$.
    From \eqref{rel}, one has
    \begin{equation*}\int_{S^{n-1}}\int_{S^{n-1}}[u,v]^p\rho_K(u)^{n+p}\rho_L(v)^{n+p}\,du\,dv=\frac{(n+p)\widetilde{c}_{n,p}V(L)\widetilde{V}_{-p}(K,\Lambda_p^{\circ}L)}{n\omega_n^2}.
        \label{long-0}
    \end{equation*}
    Together with the polar coordinate, the dual $L_{-p}$ Minkowski
    inequality \eqref{dm}, and  the $L_p$-sine Blaschke-Santal\'{o} inequality \eqref{sine}, one gets
    \begin{eqnarray}\int_{K}\int_{L}[x,y]^p\,dx\, dy
        \!\!&=&\!\!\bigg(\frac{n\omega_{n}}{n+p}\bigg)^2\int_{S^{n-1}}\int_{S^{n-1}}[u,v]^p\rho_K(u)^{n+p}\rho_L(v)^{n+p}\,du\,dv
         \nonumber\\ \!\!&=&\!\! \frac{n\widetilde{c}_{n,p}V(L)\widetilde{V}_{-p}(K,\Lambda_p^{\circ}L)}{n+p}
    \label{long-1}\\ \!\!&\geq&\!\!
        \frac{n\widetilde{c}_{n,p}
            V(L)V(K)^{\frac{n+p}{n}}V(\Lambda_p^{\circ}L)^{-\frac{p}{n}}}{n+p}
        \nonumber\\ \!\!&=&\!\! \frac{n\widetilde{c}_{n,p}[V(L)V(K)]^{\frac{n+p}{n}}}{n+p}
        \cdot \big[V(L) V(\Lambda_p^{\circ}L)\big]^{-\frac{p}{n}}\nonumber\\
        \!\!&\geq &\!\!
        \bigg(\frac{n\widetilde{c}_{n,p}[V(K)V(L)]^{\frac{n+p}{n}}}{n+p}\bigg) \cdot
        \omega_n^{-\frac{2p}{n}}  \nonumber  \\\!\!&=&\!\!\frac{ n(n-1) \omega_{n-1}\omega_{n+p-2}}{ (n+p)^2 \omega_{n+p-3}\omega_n^{1+2p/n}} \big[V(K)V(L)\big]^{\frac{n+p}{n}},\nonumber
    \end{eqnarray}
    which concludes the proof of inequality \eqref{gm}.

    Now the equality condition of \eqref{gm} immediately follows from the equality conditions of \eqref{sine} and \eqref{dm} together with the facts \eqref{cB} and \eqref{lam2}.
\end{proof}

Conversely, one can also deduce Theorem \ref{lpm} by applying
Theorem \ref{m1}. To see this, for $K\subset \Rn$ being an
$L_{n+p}$-star with positive volume, one has $\Lambda_p^{\circ}K\in
\mathscr{K}_{e}^n$. Taking $L=\Lambda_p^{\circ}K$ in \eqref{gm} yields
\begin{equation}\label{gmf}
    \int_K\int_{\Lambda_p^{\circ}K}[x,y]^p\,dx\,dy\geq
    \frac{ n(n-1) \omega_{n-1}\omega_{n+p-2}}{ (n+p)^2 \omega_{n+p-3}\omega_n^{1+2p/n}} \big[V(K)V(\Lambda_p^{\circ}K)\big]^{\frac{n+p}{n}}.
\end{equation}
On the other side, it follows from \eqref{long-1} and \eqref{rel21} that
\begin{equation}\label{gmg}\int_K\int_{\Lambda_p^{\circ}K}[x,y]^p\,dx\,dy
    =\frac{n \widetilde{c}_{n,p}V(\Lambda_p^{\circ}K)\widetilde{V}_{-p}(K,\Lambda_p^{\circ}\Lambda_p^{\circ}K)}{n+p}
    =\frac{n\widetilde{c}_{n,p} V(\Lambda_p^{\circ}K)V(K)}{n+p}.
\end{equation}
Now, inequality \eqref{sine} immediately follows from \eqref{gmf} and \eqref{gmg}.

 Just like the proof of \cite[Theorem A]{LZ}, for nonnegative
functions $f,g\in L^1(S^{n-1})$, by taking $\rho_K=f^{1/(n+p)}$ and
$\rho_L=g^{1/(n+p)}$ in Theorem \ref{m1}, one obtains the following result. 

\bt \label{m12}

Let $p\ge 1$ and $n\ge 2$. For any nonnegative functions $f,g\in
L^1(S^{n-1})$, one has
\begin{equation}\label{mm1}
    \int_{S^{n-1}}\int_{S^{n-1}}[u,v]^pf(u)g(v)\,du\,dv\geq
    \bigg(\frac{(n-1)\omega_{n-1}\omega_{n+p-2}}{n\omega_{n}\omega_{n+p-3}}\bigg)\cdot
    \|f\|_{\frac{n}{n+p}}\|g\|_{\frac{n}{n+p}}.
\end{equation}
When $n=2$ and $ \|f\|_{\frac{2}{2+p}}\|g\|_{\frac{2}{2+p}}>0$,
equality holds if and only if there exist $\phi\in
\mathrm{SL}(2)$~and real numbers $d_1,d_2>0$,  such that, for almost
all $u\in S^{1}$ (with respect to the spherical measure on $S^1$),
\begin{equation*}
    f(u)=d_1|\phi u|^{-2-p}\quad \mathrm{and} \quad g(u)=d_2|\phi
    u|^{-2-p}.\label{equ-fun-1}
\end{equation*}
When $n\ge 3$ and $ \|f\|_{\frac{n}{n+p}}\|g\|_{\frac{n}{n+p}}>0$, equality holds if and only
if there exist real numbers $d_3,d_4>0$, such that, for  almost all
$u\in S^{n-1}$ (with respect to the spherical measure on $\sphere$),
\begin{equation*}
    f(u)=d_3 \quad \mathrm{and} \quad g(u)=d_4. \label{equ-fun-2}
\end{equation*}\et
As  pointed out  by one of the referees,  Theorem \ref{m1} is
a special case of   \cite[Corollary 4.2]{DPP}. 
We also mention that Theorem  \ref{m12} can be viewed as the ``sine cousin"
of the following inequality proved by Lutwak and Zhang \cite{LZ}:
for $p\ge 1$ and continuous functions $f,g: S^{n-1}\rightarrow
(0,\infty)$, one has
\begin{equation}\label{BL}
    \int_{S^{n-1}}\int_{S^{n-1}}|u\cdot v|^pf(u)g(v)\,du\,dv\ge
    \frac{\omega_{n+p-2}}{\omega_2\omega_{n-2}\omega_{p-1}}\|f\|_{\frac{n}{n+p}}\|g\|_{\frac{n}{n+p}},
\end{equation}
with equality if and only if there exist $\phi\in \mathrm{SL}(n)$
and constants $d_1, d_2>0$,  such that, for all $u\in \sphere$,
$$f(u)=d_1|\phi u|^{-(n+p)}\ \ \mathrm{and}\ \
g(u)=d_2|\phi^{-t} u|^{-(n+p)}. $$

Note that a stronger version of inequality \eqref{BL} was
proved by Nguyen \cite{Ng}, where the continuity of $f, g$ is
replaced by integrability. Furthermore, Nguyen \cite{Ng} showed that
inequality \eqref{BL} can be used to prove the following $L_p$ moment-entropy
inequality established by Lutwak, Yang, and Zhang \cite{LYZ04}.
Suppose $p\ge 1$, $n\ge 2$, and $\lambda\in
(\frac{n}{n+p},\infty]$. Let $X,Y$ be two independent random vectors
in $\Rn$ with density functions $f,g: \Rn\rightarrow [0, \infty)$. If $X,Y$ have finite $p$th moment (that is,
$\int_{\Rn}|x|^p f(x)\,dx<\infty$, and $f$ is assumed to be
bounded if $\lambda=\infty$),  then
\begin{equation}\label{mm}\bbE(|X\cdot
    Y|^p)\ge \frac{c_0^2 (n+p)\omega_{n+p}}{n\pi \omega_{p-1}\omega_n^{1+2p/n}}\big[N_{\lambda}(X)N_{\lambda}(Y)\big]^{p/n}.\end{equation}
Here,   $$\bbE(|X\cdot Y|^p)=\int_{\Rn}\int_{\Rn} |x \cdot y|^pf(x)g(y)\,dx\,dy,$$  and the $\lambda$-R\'{e}nyi entropy power $ N_{\lambda}(X)$ of $X$   is defined by
\begin{eqnarray}
    N_{\lambda}(X)=\left\{ \begin{array}{ll}
        \big(\int_{\mathbb R^n}f(x)^{\lambda}\,dx\big)^{\frac{1}{1-\lambda}}, & \textrm{if}~\lambda\neq 1, \\
        \exp(-\int_{\mathbb R^n}f(x)\log f(x)\,dx), & \textrm{if}~\lambda= 1,\\
        \lim_{\lambda\rightarrow \infty}N_{\lambda}(X), & \textrm{if}~\lambda=\infty.
    \end{array}   \right. \notag
\end{eqnarray}
The precise value of  $c_0$ can be found in \cite{LYZ04}. Equality
in \eqref{mm} holds if and only if there exists an origin-symmetric
ellipsoid $\mathscr{E}$, such that, a.e.,
\begin{equation*}
    f(x)\!=\!\left\{\begin{array}{ll}\! b_1\cdot p_{\lambda}\big(\frac{a_1}{\rho_{\mathscr{E}}(x)}\big), & \mathrm{if}\  \lambda\in (\frac{n}{n+p},\infty), \\ \!\frac{\mathbf{1}_{a_1\mathscr{E}}(x)}{V(a_1\mathscr{E})},   &  \mathrm{if}\   \lambda=\infty,\end{array}\right. \! \mathrm{and} \ \
    g(y)\!=\!\left\{\begin{array}{ll}\! b_2\cdot p_{\lambda}\big(\frac{a_2}{\rho_{\mathscr{E}^{\circ}}(y)}\big), &  \mathrm{if}\  \lambda\in (\frac{n}{n+p},\infty), \\ \! \frac{\mathbf{1}_{a_2\mathscr{E}^{\circ}}(y)}{V(a_2\mathscr{E}^{\circ})},   & \mathrm{if}\  \lambda=\infty,\end{array}\right.  \label{che-6-6-lyz}
\end{equation*} where $a_1,  a_2 >0$ are some constants, and  $b_1, b_2$ are constants chosen to make $f$ and $g$ density functions. 
Here, the function $p_{\lambda}:(0,\infty)\rightarrow(0,\infty)$ is given by
\begin{eqnarray}
    p_{\lambda}(s)= \left\{\begin{array}{ll}
        (1+s^p)^{1/(\lambda-1)}, \ \ & \textrm{if} \ \ \lambda< 1, \\
        e^{-s^p}, \ \ & \textrm{if}\ \ \lambda=1,\\
        (1-s^p)_+^{1/(\lambda-1)}, \ \ & \textrm{if}\ \ \lambda>1,
    \end{array} \right. \notag
\end{eqnarray} where $a_+=\max\{a, 0\}$ for $a\in \R$.

Following from Theorem
\ref{m12} and along the approach in  \cite{Ng} (with slight modifications),  one can get the
 $L_p$-sine moment-entropy inequality in Theorem \ref{m2}. Note that this theorem can also be proved by Theorem \ref{lpm} and along the approach in \cite{LYZ04} (with slight modifications).

\bt \label{m2} Suppose $p\ge 1$, $n\ge 2$, and $\lambda\in
(\frac{n}{n+p},\infty]$.  Let $X,Y$ be two independent random vectors
in $\Rn$ with density functions $f,g: \Rn\rightarrow [0, \infty)$. If $X,Y$ have finite $p$th moment, then
\begin{equation}\label{inf}
    \bbE([X,Y]^p)\geq  \frac{c_0^2 (n-1)\omega_{n-1}\omega_{n+p-2}}{n \omega_{n+p-3}\omega_n^{1+2p/n} }  \Big[N_{\lambda}(X)
    N_{\lambda}(Y)\Big]^{\frac{p}{n}}.
\end{equation}
When $n=2$, equality in \eqref{inf} holds if and only if there
exists an origin-symmetric ellipsoid $\mathscr{E}$, such that, a.e.,
\begin{equation*}
    f(x)\!=\!\left\{\!\!\begin{array}{ll} b_1\cdot p_{\lambda}\big(\frac{a_1}{\rho_{\mathscr{E}}(x)}\big), & \mathrm{if}\  \lambda\in (\frac{2}{2+p},\infty), \\ \frac{\mathbf{1}_{a_1\mathscr{E}}(x)}{V(a_1\mathscr{E})},   &  \mathrm{if}\   \lambda=\infty,\end{array}\right.  \mathrm{and} \
    g(y)\!=\!\left\{\!\! \begin{array}{ll} b_2\cdot p_{\lambda}\big(\frac{a_2}{\rho_{\mathscr{E}}(y)}\big), &  \mathrm{if}\  \lambda\in (\frac{2}{2+p},\infty), \\ \frac{\mathbf{1}_{a_2\mathscr{E}}(y)}{V(a_2\mathscr{E})},   & \mathrm{if}\  \lambda=\infty, \end{array}\right.  \label{che-6-6}
\end{equation*} where $a_1,  a_2 >0$ are some constants, and  $b_1, b_2$ are constants chosen to make $f$ and $g$ density functions.
When $n\ge 3$, equality in \eqref{inf} holds if and only if, a.e.,
\begin{equation*}
    f(x)\!=\!\!\left\{\!\!\begin{array}{ll} b_3 \cdot p_{\lambda}(\frac{a_3}{\rho_{\ball}(x)}),   & \mathrm{if}\  \lambda\in (\frac{n}{n+p},\infty), \\
        \frac{\mathbf{1}_{a_3B^n}(x)}{V(a_3B^n)}, &\mathrm{if} \
        \lambda=\infty,
    \end{array}\right.   \mathrm{and} \
    g(y)\!=\!\left\{\!\!\begin{array}{ll} b_4 \cdot p_{\lambda}(\frac{a_4}{\rho_{\ball}(y)}),   & \mathrm{if}\  \lambda\in (\frac{n}{n+p},\infty),  \\
        \frac{\mathbf{1}_{a_4B^n}(y)}{V(a_4B^n)}, &\mathrm{if} \
        \lambda=\infty,
    \end{array}\right.   \label{che-6-7}
\end{equation*} where $a_3,  a_4 >0$ are some constants, and  $b_3, b_4$ are constants chosen to make $f$ and $g$ density functions.
\et

If $X$ and $Y$ are random vectors with density functions
$\mathbf{1}_K/V(K)$ and $\mathbf{1}_L/V(L)$, then inequality
\eqref{inf} reduces to inequality \eqref{gm} and Theorem \ref{m1}
can be extended to the case that $K$ and $L$ are bounded and
measurable sets in $\Rn$. In particular,  taking  $p\rightarrow
\infty$ on both sides of inequality \eqref{gm},  Stirling's formula
yields that
\begin{equation}\label{general-bs} \sup_{x\in K, y\in L}[x,y]\geq
\omega_n^{-2/n} [V(K)V(L)]^{1/n}.\end{equation} Note that the sine
Blaschke-Santal\'{o} inequality \eqref{sine-1} is a special case of
inequality \eqref{general-bs}. Indeed, under the assumptions on $K$
in Theorem \ref{santalo-sine-12-27} and $L=\ksp$, one has
\begin{equation*}   1=\sup_{x\in K, y\in \ksp}[x,y]\geq
    \omega_n^{-2/n} [V(K)V(\ksp)]^{1/n},\end{equation*} as desired.

Finally, it is easy to see that inequalities \eqref{mm1} and
\eqref{BL}, \eqref{inf} and \eqref{mm}, together with their equality
conditions are in fact equivalent when $n=2$.

\vskip 2mm \noindent  {\bf Acknowledgement.}  The authors would like
to thank Professor Daniel Hug from Karlsruhe Institute of Technology
for many helpful discussion on the proof of inequality
\eqref{compare-polar-diamond}. The authors are also indebted to the
referees for many valuable suggestions and comments, which  greatly  improve the quality and presentation of the present paper. In particular, we
are grateful to the referees for pointing out references \cite{ASW, DPP,MW-3} and providing the Mathematica code for the
images in Figure 1.

 The research of QH was supported by
NSFC (No. 11701219) and AARMS postdoctoral fellowship (joint with
Memorial University of Newfoundland, Canada). The research of AL was
supported by  Zhejiang Provincial Natural Science Foundation of
China (LY22A010001). The research of DX was supported by NSFC (No.
12071277) and  STCSM program (No. 20JC1412600). The research of DY
was supported by a NSERC grant, Canada.

\bibliographystyle{amsalpha}

\vskip 2mm \noindent
Qingzhong Huang,   \ {\small \tt hqz376560571@163.com} \\
College of Data Science, Jiaxing University, Jiaxing, 314001, China

\vskip 2mm \noindent Ai-Jun Li,  \  {\small \tt liaijun72@163.com}\\
School of Science,
 Zhejiang University of Science and Technology, Hangzhou, Zhejiang, 310023, China

 \vskip 2mm \noindent Dongmeng Xi, \  {\small \tt dongmeng.xi@live.com} \\  Department of Mathematics, Shanghai University, Shanghai, 200444, China

\vskip 2mm \noindent  Deping Ye,  \  {\small \tt deping.ye@mun.ca} \\
 Department of Mathematics and Statistics, Memorial
University of Newfoundland, St. John's, Newfoundland A1C 5S7,
Canada

\end{document}